\documentclass[12pt]{article}

\usepackage{amsmath,amssymb,amsthm,mathtools}
\usepackage{bm}
\usepackage{fullpage}

\newcommand{\R}{\mathbb R}

\newcommand{\E}{\mathbb E}

\newtheorem{theorem}{Theorem}
\newtheorem{proposition}[theorem]{Proposition}
\newtheorem{lemma}[theorem]{Lemma}
\newtheorem{corollary}[theorem]{Corollary}

\newtheorem{definition}[theorem]{Definition}

\title{Global well-posedness and regularization for a nonlinear
relativistic Fokker--Planck equation}
\author
 {
       Simone Calogero  \\
       {\small Department of Mathematical Sciences}  \\
       {\small Chalmers University of Technology}  \\
       {\small Gothenburg, Sweden} \\
       }

\date{}

\begin{document}

\maketitle

\begin{abstract}
We consider a nonlinear relativistic Fokker--Planck equation on
momentum space, whose diffusion coefficient and drift vector field
are nonlinear integral functionals of the solution.
We prove global existence and uniqueness of weak probability density
solutions preserving the relativistic
energy-momentum $(E,P)$.
A fundamental step in the proof is to establish that the scalar diffusion coefficient remains uniformly positive on every finite time interval.
The key mechanism ensuring this property is the strict relativistic
inequality $E^2-|P|^2>1$, combined with the uniform integrability of
the energy tails on every finite time interval.
We further prove that global weak solutions become
instantaneously smooth and strictly positive for positive times.
Finally, for even initial data, we prove that the diffusion coefficient remains uniformly
positive for all time.
\end{abstract}

%===========================================================
\section{Introduction}
%===========================================================

In a recent paper~\cite{Calogero}, the author introduced a nonlinear
kinetic Fokker--Planck equation which is Lorentz invariant and preserves
particle number, energy, and momentum. The purpose of the present paper
is to establish a global Cauchy theory for the spatially homogeneous
version of this equation. A fundamental analytical difficulty is that
both the drift vector field and the scalar diffusion coefficient depend
nonlocally and nonlinearly on the solution.
In particular, the diffusion coefficient is not a priori
bounded away from zero. A key step in the proof of
global existence is therefore to establish a uniform
positive lower bound for this coefficient on every finite
time interval.

The equation under consideration is
\begin{equation}\label{mainPDE}
\partial_tG
=
\Phi_G\nabla_p\cdot(D\nabla_pG)
+
\nabla_p\cdot(W_{G}G).
\end{equation}
Here, $G=G(t,p)$ is the probability density at time $t\geq0$ for a
system of particles with momentum $p\in\R^3$. The diffusion matrix
$D(p)$ is given by
\begin{equation}\label{Ddef}
D(p)
=
\frac{I+p\otimes p}{\langle p\rangle},
\qquad
\langle p\rangle=\sqrt{1+|p|^2}.
\end{equation}
The scalar diffusion coefficient $\Phi_G(t)$ and the drift vector field
$W_G(t,p)$ are given by
\begin{equation}\label{Phidef}
\Phi_G(t)
=
\frac{\mu}{3}
\left[
1-
\int_{\R^3}
\frac{G(t,p)}
{\langle p\rangle\beta_G(t,p)}\,dp
\right],\qquad
W_{G}(t,p)
=
\frac{\mu}{\langle p\rangle}
\left(
p-\frac{j_G(t)}{\beta_G(t,p)}
\right),
\end{equation}
where $\mu>0$ is the drift constant and
\begin{equation}\label{adef}
\beta_G(t,p)
=
\langle p\rangle\bigl(1-j_G(t)\cdot v(p)\bigr),
\qquad
j_G(t)=\int_{\R^3}v(p)\,G(t,p)\,dp,
\qquad
v(p)=\frac{p}{\langle p\rangle}.
\end{equation}
Thus, both the diffusion coefficient and the drift vector field are
integral nonlinear functionals of the solution $G$.
We employ units in which the speed of light and the particle mass are
equal to one.

For a probability density $G$, we define the energy $E_G$ and the momentum $P_G$ as
\begin{equation}\label{energy-mom}
E_G=\int_{\R^3}\langle p\rangle G(p)\,dp,\qquad P_G=\int_{\R^3}p\, G(p)\,dp,
\end{equation}
and the energy tail by
\begin{equation}\label{energytaildef}
\mathcal T_R(G)
=
\int_{\langle p\rangle>R}
\langle p\rangle G(p)\,dp,
\qquad R>1.
\end{equation}

Our first main result gives global well-posedness of the Cauchy problem
in the class of probability densities with finite energy and,
at the same time, a quantitative estimate excluding finite time
degeneration of the nonlinear diffusion coefficient.

\begin{theorem}[Global well-posedness]\label{mainthm}
Let $G_0\in L^1(\R^3)$, $E>1$, $P\in\R^3$ satisfy
\[
G_0\geq0\ \text{a.e.},\qquad
\int_{\R^3}G_0(p)\,dp=1,\qquad
E_{G_0}=E,\qquad
P_{G_0}=P.
\]
Then the Cauchy problem for~\eqref{mainPDE}
with initial data $G_0$ has a unique global weak solution
\[
G\in C([0,\infty);L^1(\R^3))
\]
such that
\[
G(t,\cdot)\geq0\ \text{a.e.},\qquad
\int_{\R^3}G(t,p)\,dp=1,
\qquad t\geq0.
\]
The energy and momentum are conserved:
\[
E_{G(t)}=E,
\qquad
P_{G(t)}=P,
\qquad t\geq0.
\]
Moreover,
\[
\Delta:=E^2-|P|^2-1>0.
\]
The energy tails are uniformly integrable on every finite time interval; that is, for all $T>0$,  
\begin{equation}\label{mainenergytails}
\lim_{R\to\infty}
\sup_{0\leq t\leq T}
\mathcal T_R(G(t))=0.
\end{equation}
Finally, for every $R>1$ and $t\geq0$,
\begin{equation}\label{mainPhibound}
\Phi_G(t)
\geq
\frac{9\mu}{8R^2
}
\left(
\frac{\Delta}{4E}-\mathcal T_R(G(t))
\right)_+^2,
\end{equation}
where $(x)_+=\max\{0,x\}$. 
In particular, by~\eqref{mainenergytails}, for every $T>0$,
there exists $\delta_T>0$ such that
\[
\Phi_G(t)\geq\delta_T,
\qquad 0\leq t\leq T.
\]
\end{theorem}

The main steps in the proof of Theorem~\ref{mainthm} are the following.
First, local existence is obtained by solving a fixed
point problem on the nonlinear functionals
$(\Phi_G,j_G)$.
The local solution can be continued globally
provided that, on every finite time interval,
the energy of the solution remains bounded and
$\Phi_G$ remains bounded away from zero.
The mechanism preventing the breakdown of local solutions is a
consequence of the relativistic structure of the equation.
We first establish conservation of energy $E$ and momentum $P$
directly at the level of weak solutions and prove the uniform
integrability property~\eqref{mainenergytails} of the energy tails.
The strict relativistic inequality
$E^2-|P|^2>1$ yields a quantitative lower bound
for the velocity variance in terms of the energy tails.
A coercivity estimate then gives~\eqref{mainPhibound}.
Combined with~\eqref{mainenergytails}, this estimate
shows that the diffusion coefficient remains uniformly
positive on every finite time interval, thereby allowing
the local solution to be continued globally.

The second main result shows that the global weak solutions immediately
acquire full regularity and become strictly positive.

\begin{theorem}[Regularization and strict positivity]
\label{regularizationthm}
Let $G$ be the global weak solution given by
Theorem~\ref{mainthm}. Then
\[
G\in C^\infty((0,\infty)\times\R^3).
\]
Moreover,
\[
G(t,p)>0
\qquad
\text{for every }t>0,\ p\in\R^3.
\]
\end{theorem}

The proof of Theorem~\ref{regularizationthm} involves a second
difficulty: the diffusion matrix $D(p)$ is not uniformly elliptic
on the whole momentum space. It is, however, uniformly elliptic
on every bounded momentum set. We first derive evolution identities
for the nonlinear coefficients and use a class of admissible moment
observables to prove that these coefficients are smooth for positive
times. Once this regularity has been established, the equation
satisfied by a given weak solution can be regarded as a linear,
nonautonomous parabolic equation with smooth coefficients that is
locally uniformly elliptic in momentum. Standard interior parabolic
regularity then yields instantaneous smoothness of the solution.
Finally, the parabolic strong maximum principle
gives strict positivity for every positive time.

Our third result concerns the long time behavior of the
nonlinear diffusion coefficient. For general initial data,
the estimate~\eqref{mainPhibound} does not exclude the possibility that
$\Phi_G(t)$ approaches zero as $t\to\infty$.
In the following theorem, we prove that for all even initial data
with finite energy, the diffusion coefficient remains
uniformly positive for all time.
Moreover, if the conserved energy satisfies $1<E<2$,
the energy tails are uniformly integrable on $[0,\infty)$.

\begin{theorem}[Uniform positivity of the diffusion coefficient for even initial data]
\label{uniformPhitheorem}
Let $G_0\in L^1(\R^3)$ satisfy
\[
G_0\geq0\ \text{a.e.},\qquad
G_0(p)=G_0(-p)\ \text{a.e.},\qquad
\int_{\R^3}G_0(p)\,dp=1,
\]
and assume that $E_{G_0}=E<\infty$. 
Let $G$ be the global weak solution given by
Theorem~\ref{mainthm}. Then
\[
G(t,p)=G(t,-p),
\qquad t>0,
\]
and there exists $\delta>0$, depending only
on $G_0$ and $\mu$, such that
\begin{equation}\label{uniformPhibound}
\Phi_G(t)\geq\delta,
\qquad t\geq0.
\end{equation}
If, in addition, $E<2$, then
\begin{equation}\label{globaluniformtails}
\lim_{R\to\infty}
\sup_{t\geq0}\mathcal T_R(G(t))=0.
\end{equation}
\end{theorem}

The key observation in the proof of
Theorem~\ref{uniformPhitheorem} is that the number
current $j_G$ vanishes identically for even solutions,
in which case the particle energy distribution
satisfies a one-dimensional Fokker--Planck equation.
When $1<E<2$, a PDE comparison argument with an
explicit stationary supersolution yields the uniform
integrability property~\eqref{globaluniformtails}.
For arbitrary finite energy, we establish a comparison
principle for the particle energy distributions of
even solutions. This allows us to compare the diffusion
coefficient with that of a radially nonincreasing solution
whose conserved energy is less than two, thereby proving
the uniform lower bound~\eqref{uniformPhibound}.

The validity of the lower bound
\[
\inf_{t\geq0}\Phi_G(t)>0
\]
for general initial data remains an open problem, as is the
convergence of solutions to equilibrium.
The latter problem presents additional difficulties compared
with the linear case studied in~\cite{Calogero3,Calogero2},
in part because the equilibrium distribution
of~\eqref{mainPDE} decays only polynomially.
Indeed, as shown in~\cite{Calogero}, the equilibrium
distribution of~\eqref{mainPDE} is not the standard exponentially
decaying Maxwell--J\"uttner distribution; rather, it belongs
to the family
\[
G_{\mathrm{eq}}(p)
=
\frac{Z(\lambda,y)}
{(\langle p\rangle-y\cdot p)^{2\lambda}},
\qquad
\lambda>2,\quad |y|<1,
\]
where
\[
Z(\lambda,y)
=
\frac{
(1-|y|^2)^{\lambda+1/2}\Gamma(\lambda)
}{
\pi^{3/2}\Gamma(\lambda-3/2)
}
\]
is the normalization constant.
For these distributions,
\[
j_{G_{\mathrm{eq}}}=y,
\qquad
\Phi_{G_{\mathrm{eq}}}=\frac{\mu}{2\lambda}.
\]
The parameters $\lambda$ and $y$ are related to the
conserved energy and momentum of the solution.
In particular, when $y=0$, the equilibrium is isotropic
and takes the form
\[
G_{\mathrm{eq}}(p)
=
\frac{\Gamma(\lambda)}
{\pi^{3/2}\Gamma(\lambda-3/2)}
\langle p\rangle^{-2\lambda}.
\]
The full asymptotic behavior of solutions to~\eqref{mainPDE}
will be addressed in a subsequent publication.

%===========================================================
\section{Preliminaries}
%===========================================================

We begin by noticing that, for $|p|\leq R$ and
$\xi\in\R^3$, the diffusion matrix~\eqref{Ddef} satisfies
\begin{equation}\label{localellipticity}
\frac{1}{\sqrt{1+R^2}}|\xi|^2
\leq
\xi^TD(p)\xi
\leq
\sqrt{1+R^2}|\xi|^2.
\end{equation}
Thus, the operator
$G\to\nabla_p\cdot(D\nabla_pG)$ in the right hand
side of~\eqref{mainPDE} is uniformly elliptic on
every bounded momentum set.

\begin{definition}[Weak probability density solution]\label{weaksolutiondefinition}
Let $T>0$. A function
\[
G\in C([0,T];L^1(\R^3))
\]
is a weak probability density solution of~\eqref{mainPDE} on $[0,T]$
with initial data $G_0$ if
\[
G(t,\cdot)\geq0\ \text{a.e.},
\qquad
\int_{\R^3}G(t,p)\,dp=1,
\qquad 0\leq t\leq T,
\]
and, for every
$\psi\in C_c^\infty([0,T)\times\R^3)$,
\begin{equation}\label{weakformulation}
\int_0^T\int_{\R^3}
G\left[
\partial_t\psi
+\Phi_G\nabla_p\cdot(D\nabla_p\psi)
-W_G\cdot\nabla_p\psi
\right]dp\,dt
+
\int_{\R^3}G_0(p)\psi(0,p)\,dp
=0.
\end{equation}
\end{definition}

The following lemma will be used to extend local weak
solutions to larger time intervals.

\begin{lemma}\label{integratedweaklemma}
Let $G\in C([0,T];L^1(\R^3))$ be a weak
solution of~\eqref{mainPDE} in the sense of
Definition~\ref{weaksolutiondefinition}.
Then, for every $0\leq t_1<t_2\leq T$
and every
$\psi\in C^{1,2}([0,T]\times\R^3)$
with compact support in the momentum variable
uniformly on $[0,T]$, there holds
\begin{align}
&\int_{\R^3}G(t_2,p)\psi(t_2,p)\,dp
-\int_{\R^3}G(t_1,p)\psi(t_1,p)\,dp
\nonumber\\
&\quad=
\int_{t_1}^{t_2}\int_{\R^3}
G(t,p)
\left[
\partial_t\psi
+\Phi_G\nabla_p\cdot(D\nabla_p\psi)
-W_G\cdot\nabla_p\psi
\right]dp\,dt.
\label{integratedweak}
\end{align}
In particular, the restriction of a weak solution
to any subinterval $[t_0,T]$ is a weak solution
with initial data $G(t_0,p)$ almost everywhere.
Moreover, if $G_1$ and $G_2$ are weak solutions
on $[0,t_0]$ and $[t_0,T]$, respectively, and
$G_1(t_0,p)=G_2(t_0,p)$ almost everywhere,
then their concatenation is a weak solution
on $[0,T]$.
\end{lemma}

\begin{proof}
Let $\psi$ be as in the statement and define
\[
F(t)=
\int_{\R^3}G(t,p)\psi(t,p)\,dp
\]
and
\[
B(t)=
\int_{\R^3}G(t,p)
\left[
\partial_t\psi
+\Phi_G\nabla_p\cdot(D\nabla_p\psi)
-W_G\cdot\nabla_p\psi
\right]dp.
\]

Since $G\in C([0,T];L^1(\R^3))$ and $|v(p)|<1$,
the function $j_G$ is continuous and satisfies
\[
|j_G(t)|
\leq
\int_{\R^3}|v(p)|G(t,p)\,dp<1,
\qquad 0\leq t\leq T.
\]
By compactness, there exists $\gamma_T<1$ such that
\[
|j_G(t)|\leq\gamma_T,
\qquad 0\leq t\leq T.
\]
It follows directly from the definitions of
$\Phi_G$ and $W_G$ in~\eqref{Phidef} that these coefficients
are continuous in time, uniformly on compact
momentum sets. Consequently, both $F$ and $B$
are continuous on $[0,T]$.

Let $\eta\in C_c^\infty([0,T))$.
Using $\eta(t)\psi(t,p)$ in
\eqref{weakformulation}, after approximation
by smooth test functions, gives
\[
\int_0^T\eta'(t)F(t)\,dt
+
\int_0^T\eta(t)B(t)\,dt
+
\eta(0)I_0
=0,
\]
where
\[
I_0=\int_{\R^3}G_0(p)\psi(0,p)\,dp.
\]
First take $\eta\in C_c^\infty((0,T))$.
Then $F'=B$ in the sense of distributions
on $(0,T)$.
Since $F$ and $B$ are continuous on $[0,T]$,
it follows that
\[
F(t)=F(0)+\int_0^tB(s)\,ds,
\qquad 0\leq t\leq T.
\]
Substituting this identity into the preceding
weak identity and integrating by parts gives
\[
\eta(0)\bigl(I_0-F(0)\bigr)=0.
\]
Choosing $\eta$ with $\eta(0)=1$, we obtain
\[
\int_{\R^3}G(0,p)\psi(0,p)\,dp
=
\int_{\R^3}G_0(p)\psi(0,p)\,dp.
\]
In particular, taking $\psi(t,p)=\varphi(p)$
for arbitrary $\varphi\in C_c^\infty(\R^3)$
shows that
\[
G(0,p)=G_0(p)
\qquad\text{for almost every }p\in\R^3.
\]
Finally, subtracting the identity
\[
F(t)=F(0)+\int_0^tB(s)\,ds
\]
at $t_1$ and $t_2$ gives
\eqref{integratedweak}.
Taking $t_1=t_0$ shows that the restriction
of $G$ to $[t_0,T]$ satisfies the integrated
weak formulation with initial data $G(t_0,p)$.
In particular, it satisfies the corresponding
weak formulation after translating the initial
time to $t_0$.

Finally, suppose that $G_1$ and $G_2$ are
weak solutions on $[0,t_0]$ and $[t_0,T]$
with $G_1(t_0)=G_2(t_0)$.
Define
\[
G(t,p)=
\begin{cases}
G_1(t,p),&0\leq t\leq t_0,\\
G_2(t,p),&t_0\leq t\leq T.
\end{cases}
\]
Then $G\in C([0,T];L^1(\R^3))$.
Applying~\eqref{integratedweak} on
$[0,t_0]$ and $[t_0,T]$ and adding the
resulting identities, the boundary terms
at $t_0$ cancel.
Hence, $G$ satisfies~\eqref{integratedweak}
on $[0,T]$. Taking $t_1=0$ and integrating
by parts in time recovers the weak formulation~\eqref{weakformulation}
in Definition~\ref{weaksolutiondefinition}.
Therefore, $G$ is a weak solution on $[0,T]$.
\end{proof}

For a probability density $G$, we define
\[
j_G=\int_{\R^3}v(p)G(p)\,dp,
\qquad
\beta_G(p)=\langle p\rangle
\bigl(1-j_G\cdot v(p)\bigr),
\]
and
\[
\Phi_G
=
\frac{\mu}{3}
\left[
1-
\int_{\R^3}
\frac{G(p)}
{\langle p\rangle\beta_G(p)}\,dp
\right].
\]
Recall also the definition~\eqref{energy-mom} of the energy-momentum $(E_G,P_G)$. 
If $G$ has finite energy, then
\[
E_G>1,
\]
since $\langle p\rangle\geq1$, with equality
only at $p=0$, and $G$ is an $L^1$
probability density. 
For $E>1$, we set
\[
\gamma_E=\sqrt{1-E^{-2}}.
\]
\begin{lemma}\label{estimatesBeta}
Let $G\in L^1(\R^3)$ be a probability density. Then
\begin{equation}\label{betanoen}
|j_G|<1,\quad \frac{1}{2\langle p\rangle}\leq\beta_G(p)\leq2\langle p\rangle.
\end{equation}
If, in addition, $G$ has finite energy, then 
\begin{equation}\label{jenergybound}
|j_G|\leq\gamma_{E_G}.
\end{equation}
\end{lemma}

\begin{proof}
Since $|v(p)|<1$ for all $p\in\R^3$ and
$G$ is a probability density, we have
\[
|j_G|
\leq\int_{\R^3}|v(p)|G(p)\,dp
<\int_{\R^3}G(p)\,dp=1.
\]
The upper bound on $\beta_G(p)$  in~\eqref{betanoen} then follows immediately. For the lower bound on $\beta_G(p)$ we write
\[
\frac1{\langle p\rangle\beta_G(p)}
=
\frac{1-|v(p)|^2}{1-j_G\cdot v(p)}
\leq
\frac{1-|v(p)|^2}{1-|v(p)|}
=
1+|v(p)|
\leq2.
\]

Now assume that $G$ has finite energy.
By Jensen's inequality,
\[
|j_G|^2
=
\left|\int_{\R^3}v(p)G(p)\,dp\right|^2
\leq
\int_{\R^3}|v(p)|^2G(p)\,dp.
\]
Since $1-|v(p)|^2=\langle p\rangle^{-2}$,
we obtain
\begin{equation}\label{jGtemp}
1-|j_G|^2
\geq
\int_{\R^3}\langle p\rangle^{-2}G(p)\,dp,
\end{equation}
Applying Jensen's inequality to the convex
function $y\mapsto y^{-2}$ gives
\[
\int_{\R^3}\langle p\rangle^{-2}G(p)\,dp
\geq
\left(
\int_{\R^3}\langle p\rangle G(p)\,dp
\right)^{-2}
=
\frac1{E_G^2}.
\]
Consequently,~\eqref{jGtemp} gives
\[
1-|j_G|^2\geq\frac1{E_G^2},
\]
which proves~\eqref{jenergybound}.
\end{proof}

Bounds on the diffusion coefficient $\Phi_G$ are particularly important, so we state them in a separate lemma. 

\begin{lemma}\label{Phiboundslemma}
Let $G\in L^1(\R^3)$ be a probability density. Then
\begin{equation}\label{Phinoen}
0<\frac{9\mu}{32}
\int_{\R^3}|v(p)-j_G|^2G(p)\,dp\leq \Phi_G\leq\frac{\mu}{3}.
\end{equation}
Moreover, if $G$ has finite energy, then
\begin{equation}\label{elementaryPhi}
\Phi_G
\leq
\frac{\mu}{3}\,\gamma_{E_G}^2.
\end{equation}
\end{lemma}

\begin{proof}
The upper bound in~\eqref{Phinoen} is immediate, because $\beta_G(p)>0$. For the lower bound, set $j=j_G$. As in the definition of $\Phi_G$,
\[
\frac{3\Phi_G}{\mu}
=
\int_{\R^3}
\frac{|v|^2-j\cdot v}{1-j\cdot v}
G\,dp.
\]
Using
\[
|v|^2-j\cdot v
=
|v-j|^2+j\cdot(v-j)
\]
and
\[
1-j\cdot v
=
1-|j|^2-j\cdot(v-j),
\]
we obtain
\begin{equation}\label{identityPhi}
\frac{3\Phi_G}{\mu}
=
\int_{\R^3}
\frac{
|v-j|^2+
\dfrac{[j\cdot(v-j)]^2}{1-|j|^2}
}{1-j\cdot v}
G\,dp.
\end{equation}
We claim that
\begin{equation}\label{jineq}
\frac{
|v-j|^2+
\dfrac{[j\cdot(v-j)]^2}{1-|j|^2}
}{1-j\cdot v}
\geq
\frac{|v-j|^2}{(1-|j|)(1+|j|)^2}.
\end{equation}
For $j=0$, or $v=j$, this is an equality. Otherwise,
set
\[
a=|j|,\qquad
u=\frac{j\cdot v}{a},
\qquad r=|v|.
\]
Then $|u|\leq r<1$, and
\begin{equation}\label{temporalino}
\frac{
|v-j|^2+
\dfrac{[j\cdot(v-j)]^2}{1-a^2}
}{(1-j\cdot v)|v-j|^2}=
\frac1{1-au}
\left(
1+
\frac{a^2(u-a)^2}
{(1-a^2)(r^2+a^2-2au)}
\right).
\end{equation}
Since $r<1$, we have
\[
r^2+a^2-2au\leq1+a^2-2au.
\]
Using the identity
\[
(1-a^2)(1+a^2-2au)+a^2(u-a)^2
=(1-au)^2,
\]
we deduce that the right hand side
of~\eqref{temporalino} is bounded below by
\[
\frac{1-au}
{(1-a^2)(1+a^2-2au)}.
\]
The last expression is increasing in $u$,
since
\[
\frac{d}{du}
\left(
\frac{1-au}
{(1-a^2)(1+a^2-2au)}
\right)
=
\frac{a}{(1+a^2-2au)^2}>0.
\]
Its minimum for $-1\leq u\leq1$ is therefore
attained at $u=-1$, where it equals
\[
\frac1{(1-a)(1+a)^2}.
\]
This proves the pointwise inequality~\eqref{jineq}. 
Integrating against $G$ gives
\[
\Phi_G
\geq
\frac{\mu}{3(1-|j_G|)(1+|j_G|)^2}
\int_{\R^3}|v-j_G|^2G\,dp.
\]
Since the function
\[
a\longmapsto(1-a)(1+a)^2
\]
is increasing on $[0,1/3]$ and decreasing
on $[1/3,1)$, we obtain
\[
\Phi_G
\geq\frac{9\mu}{32}
\int_{\R^3}|v(p)-j_G|^2G(p)\,dp.
\]
The integral is strictly positive because
$v$ is injective and $G$ is an $L^1$
probability density. Hence, the lower bound in~\eqref{Phinoen} follows. 

Next assume that $G$ has finite energy. 
To prove the upper bound on $\Phi_G$ in this case, we write, using the Cauchy--Schwarz
inequality, 
\[
\left(\int_{\R^3}\frac{G}{\langle p\rangle}\,dp\right)^2
\leq
\left(\int_{\R^3}
\frac{\beta_G}{\langle p\rangle}G\,dp\right)\left(\int_{\R^3}
\frac{G}{\langle p\rangle\beta_G}\,dp\right)=
(1-|j_G|^2)
\int_{\R^3}
\frac{G}{\langle p\rangle\beta_G}\,dp.
\]
Moreover, Jensen's inequality implies
\[
\int_{\R^3}\frac{G}{\langle p\rangle}\,dp
\geq\frac1{E_G}.
\]
Hence
\[
\int_{\R^3}
\frac{G}{\langle p\rangle\beta_G}\,dp
\geq
\frac1{E_G^2(1-|j_G|^2)}
\geq\frac1{E_G^2}.
\]
Using this inequality in the definition of $\Phi_G$ gives~\eqref{elementaryPhi}.  
\end{proof}

{\bf Remark.} The upper bound in~\eqref{elementaryPhi} is sharp
for every fixed energy $E_G>1$.
Indeed, it can be approached by spherically symmetric
probability densities concentrated near the sphere
$\langle p\rangle=E_G$.

\begin{lemma}\label{Minkowskilemma}
Every probability density $G$ with finite
first energy satisfies
\begin{equation}\label{Minkowskiineq}
E_G^2-|P_G|^2>1.
\end{equation}
\end{lemma}

\begin{proof}
The function
\[
p\longmapsto\langle p\rangle=\sqrt{1+|p|^2}
\]
is strictly convex. Jensen's inequality gives
\[
E_G
\geq
\sqrt{
1+
\left|
\int_{\R^3}p\,G(p)\,dp
\right|^2
}
=
\sqrt{1+|P_G|^2}.
\]
Equality in Jensen's inequality would imply
that $p$ is constant almost everywhere
with respect to the probability distribution,
which would make the corresponding probability
measure a Dirac mass. Since $G$ is an $L^1$
probability density, this is impossible.
Therefore the inequality is strict.
\end{proof}

We continue with a simple calculus lemma.

\begin{lemma}\label{velocitylemma}
Let $v,w\in\R^3$ satisfy $|v|<1$ and $|w|<1$.
Then
\begin{equation}\label{velocityineq}
0\leq
1-v\cdot w
-\sqrt{(1-|v|^2)(1-|w|^2)}
\leq
|v-w|.
\end{equation}
\end{lemma}

\begin{proof}
Set
\[
r=|v|,\qquad s=|w|.
\]
Choose $z\in[-1,1]$ such that
\[
v\cdot w=rsz.
\]
When $rs=0$, any $z\in[-1,1]$ may be chosen.
Set
\[
A=
1-rsz-\sqrt{(1-r^2)(1-s^2)}.
\]
The first inequality in~\eqref{velocityineq}
follows from
\[
1-rsz\geq1-rs
\geq\sqrt{(1-r^2)(1-s^2)}.
\]
To prove $A\leq|v-w|$ we calculate
\[
|v-w|^2-A^2
=
r^2s^2(1-z^2)
+
2\sqrt{(1-r^2)(1-s^2)}
\left[
1-rsz-\sqrt{(1-r^2)(1-s^2)}
\right].
\]
Both terms on the right hand side are
nonnegative, hence $A^2\leq|v-w|^2$,
which proves the second inequality
in~\eqref{velocityineq}.
\end{proof}

Recall that the energy tail of a probability density $G$ is defined by~\eqref{energytaildef}.

\begin{proposition}\label{variancegap}
Let $E>1$ and $P\in\R^3$ be given such that
\[
\Delta:=E^2-|P|^2-1>0.
\]
Then every probability density $G$ satisfying $(E_{G},P_G)=(E,P)$ 
satisfies, for every $R>1$,
\begin{equation}\label{variancegapineq}
\int_{\R^3}|v(p)-j_G|^2G(p)\,dp
\geq
\frac{
\bigl(\Delta-2E(1+\gamma_E)\mathcal T_R(G)\bigr)_+^2
}{
4E^2R^2
}.
\end{equation}
\end{proposition}

\begin{proof}
For $p,q\in\R^3$, write
\[
x=\langle p\rangle,\qquad
y=\langle q\rangle,\qquad
v=v(p),\qquad w=v(q).
\]
Since $p=xv$, we have
\[
E^2-|P|^2
=
\iint_{\R^3\times\R^3}
xy(1-v\cdot w)
G(p)G(q)\,dp\,dq.
\]
Moreover,
\[
1-|v|^2=\frac1{x^2},
\qquad
1-|w|^2=\frac1{y^2},
\]
and therefore
\[
xy\sqrt{(1-|v|^2)(1-|w|^2)}=1.
\]
Since $G$ is a probability density,
\[
1=
\iint_{\R^3\times\R^3}
G(p)G(q)\,dp\,dq.
\]
Consequently,
\[
\Delta
=
\iint_{\R^3\times\R^3}
xy\left[
1-v\cdot w
-\sqrt{(1-|v|^2)(1-|w|^2)}
\right]
G(p)G(q)\,dp\,dq.
\]
By Lemma~\ref{velocitylemma},
\[
\Delta
\leq
\iint_{\R^3\times\R^3}
xy|v-w|G(p)G(q)\,dp\,dq.
\]
Using
\[
|v-w|
\leq|v-j_G|+|w-j_G|,
\]
we obtain
\[
\Delta
\leq
2E\int_{\R^3}
\langle p\rangle|v(p)-j_G|G(p)\,dp.
\]
Set
\[
V_G=
\int_{\R^3}|v(p)-j_G|^2G(p)\,dp.
\]
Splitting the last integral into
$\langle p\rangle\leq R$ and
$\langle p\rangle>R$, and using
$|v-j_G|\leq1+\gamma_E$, which follows
from Lemma~\ref{estimatesBeta}, gives
\[
\begin{aligned}
\int_{\R^3}
\langle p\rangle|v-j_G|G\,dp
&\leq
R\int_{\R^3}|v-j_G|G\,dp
+
(1+\gamma_E)\mathcal T_R(G)\\
&\leq
R\sqrt{V_G}
+
(1+\gamma_E)\mathcal T_R(G).
\end{aligned}
\]
Hence
\[
\Delta
\leq
2ER\sqrt{V_G}
+
2E(1+\gamma_E)\mathcal T_R(G).
\]
Rearranging yields
\[
\sqrt{V_G}
\geq
\frac{
\bigl(\Delta-2E(1+\gamma_E)
\mathcal T_R(G)\bigr)_+
}{
2ER
},
\]
which proves~\eqref{variancegapineq}.
\end{proof}

\begin{corollary}\label{uniformtailPhigap}
Let $G(t,\cdot)$, $t\geq0$, be a family of
probability densities satisfying
\[
E_{G(t)}=E,
\qquad
P_{G(t)}=P,
\]
and $\Delta=E^2-|P|^2-1>0$. 
Assume that there exists $R_*>1$ such that
\begin{equation}\label{uniformtailcondition}
\sup_{t\geq0}\mathcal T_{R_*}(G(t))
<
\frac{\Delta}{2E(1+\gamma_E)}.
\end{equation}
Then
\begin{equation}\label{uniformtailPhibound}
\inf_{t\geq0}\Phi_G(t)>0.
\end{equation}
In particular,~\eqref{uniformtailPhibound} holds if
\[
\lim_{R\to\infty}
\sup_{t\geq0}\mathcal T_R(G(t))=0.
\]
\end{corollary}

\begin{proof}
By Proposition~\ref{variancegap}, the inequality~\eqref{uniformtailcondition}
implies
\[
\inf_{t\geq0}
\int_{\R^3}|v-j_G(t)|^2G(t,p)\,dp>0.
\]
The conclusion~\eqref{uniformtailPhibound}
then follows from Lemma~\ref{Phiboundslemma}.
\end{proof}

To conclude this section we prove a Lipschitz
continuity result for the coefficient functionals
$j_G$ and $\Phi_G$.
We first recall the definition of the
$1$-Wasserstein distance~\cite{Villani2003}.
Let $G$ and $H$ be probability densities on
$\R^3$ with finite energy. Then
\[
W_1(G,H)
:=
\inf_{\pi\in\Pi(G,H)}
\int_{\R^3\times\R^3}|p-q|\,d\pi(p,q),
\]
where $\Pi(G,H)$ denotes the set of couplings
of the measures
$\nu_G=G(p)\,dp$ and $\nu_H=H(p)\,dp$,
i.e., the set of probability measures on
$\R^3\times\R^3$ having $\nu_G$ and $\nu_H$
as first and second marginals, respectively. 
By the Kantorovich--Rubinstein duality theorem,
\[
W_1(G,H)
=
\sup_{\operatorname{Lip}(f)\leq1}
\left|
\int_{\R^3}f(p)(G(p)-H(p))\,dp
\right|.
\]
\begin{lemma}\label{momentLipschitz}
For every $0\leq\gamma<1$ there exists a
constant $C_\gamma$, depending also on
the fixed parameter $\mu$, such that,
whenever $G$ and $H$ are probability densities
with finite energy satisfying
\[
|j_G|,|j_H|\leq\gamma,
\]
one has
\begin{equation}\label{momentLip}
|j_G-j_H|
+
|\Phi_G-\Phi_H|
\leq
C_\gamma W_1(G,H).
\end{equation}
\end{lemma}

\begin{proof}
A direct calculation gives
\[
\partial_{p_j}v_i
=
\frac{\delta_{ij}-v_iv_j}{\langle p\rangle},
\]
and hence
\[
\nabla v(p)
=
\frac1{\langle p\rangle}
\left(I-v(p)\otimes v(p)\right).
\]
The eigenvalues of $\nabla v(p)$ are
$\langle p\rangle^{-1}$,
$\langle p\rangle^{-1}$,
$\langle p\rangle^{-3}$.
Therefore, denoting by $\|\cdot\|$
the operator norm,
\[
\|\nabla v(p)\|
=
\frac1{\langle p\rangle}
\leq1.
\]
For every $e\in\R^3$ with $|e|=1$,
the function $f(p)=e\cdot v(p)$ is
therefore $1$-Lipschitz.
By the Kantorovich--Rubinstein formula,
\[
\left|
\int_{\R^3}
f(p)(G(p)-H(p))\,dp
\right|
=
|e\cdot(j_G-j_H)|
\leq W_1(G,H).
\]
Taking the supremum over $|e|=1$ gives
\begin{equation}\label{jLip}
|j_G-j_H|\leq W_1(G,H).
\end{equation}
Set
\[
F(p,a)
=
\frac1{\langle p\rangle^2
(1-a\cdot v(p))}.
\]
Then
\[
\Phi_G-\Phi_H
=
\frac{\mu}{3}
\int_{\R^3}
\bigl(F(p,j_H)H-F(p,j_G)G\bigr)\,dp.
\]
For $|a|\leq\gamma$, there holds $1-a\cdot v\geq1-\gamma$, 
and direct differentiation gives
\[
\sup_{|a|\leq\gamma}
\|\nabla_pF(\cdot,a)\|_\infty
+
\sup_{|a|\leq\gamma}
\|\nabla_aF(\cdot,a)\|_\infty
\leq C_\gamma.
\]
Hence
\[
\begin{aligned}
\left|
\int_{\R^3}F(p,j_G)G\,dp
-
\int_{\R^3}F(p,j_H)H\,dp
\right|&\leq
\left|
\int_{\R^3}F(p,j_G)(G-H)\,dp
\right|\\
&\qquad+
\int_{\R^3}
|F(p,j_G)-F(p,j_H)|H\,dp\\
&\leq
C_\gamma W_1(G,H)
+
C_\gamma|j_G-j_H|.
\end{aligned}
\]
Together with~\eqref{jLip}, this
proves~\eqref{momentLip}.
\end{proof}

{\bf Remark.}
The coefficient functionals are also locally
Lipschitz with respect to the $L^1$ distance.
Indeed, if $|j_G|,|j_H|\leq\gamma<1$, then
\[
|j_G-j_H|+|\Phi_G-\Phi_H|
\leq C_\gamma\|G-H\|_{L^1}.
\]
The use of the $1$-Wasserstein distance in the
preceding lemma is nevertheless important for
the fixed point argument in the proof of
local existence.

%===========================================================
\section{The frozen Fokker--Planck equation}
\label{frozensec}
%===========================================================

In this section we study the linear equation obtained by freezing the
nonlinear coefficients $\Phi_G$ and $j_G$. Let
\[
\phi\in C([0,\infty)),
\qquad
a\in C([0,\infty);\R^3)
\]
be such that, for every $T>0$, there exist constants
$\delta_T,\Lambda_T,\gamma_T$ satisfying
\begin{equation}\label{frozenbounds}
0<\delta_T\leq\phi(t)\leq\Lambda_T,
\qquad
|a(t)|\leq\gamma_T<1,
\qquad 0\leq t\leq T.
\end{equation}
We consider the linear Fokker--Planck equation
\begin{equation}\label{frozenFP}
\partial_tG
=
\phi\nabla_p\cdot(D\nabla_pG)
+
\nabla_p\cdot(V_aG),
\end{equation}
where
\[
V_a(p)
=
\frac{\mu}{\langle p\rangle}
\left(
p-\frac{a}{\beta_a(p)}
\right),
\qquad
\beta_a(p)
=
\langle p\rangle
\bigl(1-a\cdot v(p)\bigr).
\]
The equation admits a useful stochastic representation. Since
\[
\sum_{j=1}^3\partial_{p_j}D_{ij}(p)=3v_i(p),
\]
 we may rewrite
\eqref{frozenFP} as
\[
\partial_tG
=
-\nabla_p\cdot(Q_{\phi,a}G)
+
\phi\,\sum_{i,j=1}^3\partial_{p_i}\partial_{p_j}(D_{ij}G),
\]
where
\begin{equation}\label{frozenB}
Q_{\phi,a}(p)
=
3\phi\,v(p)-V_a(p).
\end{equation}
Consider now the stochastic differential equation 
\begin{equation}\label{frozenSDE}
dP_t
=
Q_{\phi(t),a(t)}(P_t)\,dt
+
\sigma_{\phi(t)}(P_t)\,dB_t,
\end{equation}
where $B_t$ is a standard three dimensional Brownian motion and
\[
\sigma_\phi(p)
=
\sqrt{2\phi}\,D^{1/2}(p).
\]
An explicit positive square root of $D$ is
\[
D^{1/2}(p)
=
\langle p\rangle^{-1/2}
\left[
I+\frac{p\otimes p}{\langle p\rangle+1}
\right].
\]
The stochastic differential
equation~\eqref{frozenSDE} and the Fokker--Planck
equation~\eqref{frozenFP} have the following (standard) relation.

Suppose that $P_t$ is a solution of~\eqref{frozenSDE} with initial data $P_0$ having probability
distribution $\nu_0=G_0(p)\,dp$, and denote by $\nu_t$ the probability
distribution of $P_t$. The infinitesimal generator associated
with~\eqref{frozenSDE} is
\begin{equation}\label{GenSDE}
L_t\psi(p)
=
Q_{\phi(t),a(t)}(p)\cdot\nabla_p\psi(p)
+
\phi(t)D(p):\nabla_p^2\psi(p),
\end{equation}
where
\[
D:\nabla_p^2\psi
=
\sum_{i,j=1}^3D_{ij}\partial_{p_i}\partial_{p_j}\psi.
\]
By~\eqref{frozenB} and
\[
\nabla_p\cdot(D\nabla_p\psi)
=
D:\nabla_p^2\psi
+
3v\cdot\nabla_p\psi,
\]
we have
\begin{equation}\label{frozenGenerator}
L_t\psi
=
\phi(t)\nabla_p\cdot(D\nabla_p\psi)
-
V_{a(t)}\cdot\nabla_p\psi.
\end{equation}

Let $T>0$ and
$\psi\in C_c^\infty([0,T)\times\R^3)$.
Applying It\^o's formula to $\psi(t,P_t)$ gives
\[
\psi(T,P_T)
=
\psi(0,P_0)
+
\int_0^T
\left[
\partial_t\psi(t,P_t)
+
L_t\psi(t,P_t)
\right]dt+
\int_0^T
\nabla_p\psi(t,P_t)\cdot
\sigma_{\phi(t)}(P_t)\,dB_t.
\]
Since $\psi(T,\cdot)=0$, taking expectations yields
\begin{equation}\label{temporale}
\int_{\R^3}\psi(0,p)d\nu_0(p)
+
\int_0^T\int_{\R^3}
\left[
\partial_t\psi(t,p)
+
L_t\psi(t,p)
\right]d\nu_t(p)\,dt
=0.
\end{equation}
The stochastic integral has zero expectation because
$\nabla_p\psi$ has compact support and
$\sigma_{\phi(t)}$ is bounded on compact subsets of
$[0,T]\times\R^3$. 
Consequently, the probability measures $\nu_t$ solve the
forward Kolmogorov equation
\begin{equation}\label{forwardKol}
\partial_t\nu_t=L_t^*\nu_t
\end{equation}
in the distributional sense. 
If $\nu_t$ admits a density $G(t,p)$, then, using
\eqref{frozenGenerator}, we obtain
\begin{equation}\label{frozenWeakForm}
\int_0^T\int_{\R^3}
G(t,p)
\left[
\partial_t\psi
+
\phi(t)\nabla_p\cdot(D\nabla_p\psi)
-
V_{a(t)}\cdot\nabla_p\psi
\right]dp\,dt+
\int_{\R^3}G_0(p)\psi(0,p)\,dp
=0.
\end{equation}
This is precisely the weak formulation of the frozen
Fokker--Planck equation~\eqref{frozenFP}.
Thus, whenever the time marginal distributions of a solution
to~\eqref{frozenSDE} admit densities, these densities solve
\eqref{frozenFP} in the weak sense.

Next we show that the coefficients of
\eqref{frozenSDE} satisfy the usual bounds in the Cauchy theory for stochastic differential equations~\cite{Oksendal2003}. For a matrix $A=(A_{ij})$, we denote by
\[
\|A\|_{\rm HS}
=
\left(\sum_{i,j}|A_{ij}|^2\right)^{1/2}
\]
its Hilbert--Schmidt norm.

\begin{lemma}\label{frozencoefflemma}
Let
\[
0<\delta\leq\phi,\psi\leq\Lambda,
\qquad
|a|,|b|\leq\gamma<1.
\]
Then
\begin{equation}\label{Bbound}
|Q_{\phi,a}(p)|\leq C,\qquad \|\sigma_\phi(p)\|_{\rm HS}\leq C\langle p\rangle^{1/2}
\end{equation}
and
\begin{align}
&|Q_{\phi,a}(p)-Q_{\psi,b}(q)|
\leq
C\left(
|p-q|+|\phi-\psi|+|a-b|
\right),\label{Bstability}\\
&\|\sigma_\phi(p)-\sigma_\psi(q)\|_{\rm HS}^2
\leq
C\left(
|p-q|^2+|\phi-\psi|^2\langle q\rangle
\right).\label{sigmastability}
\end{align}
The constant $C$ depends only on
$\delta,\Lambda,\gamma$ and $\mu$. 
Finally, the coefficients $Q_{\phi,a}$ and
$\sigma_\phi$ are smooth in $p$, and for every
integer $m\geq1$ there exists a constant $C_m$,
depending only on $m,\delta,\Lambda,\gamma$ and
$\mu$, such that
\begin{equation}\label{frozenhigherderivatives}
\sup_{p\in\R^3}
\left(
|\nabla_p^m Q_{\phi,a}(p)|
+
\|\nabla_p^m\sigma_\phi(p)\|
\right)
\leq C_m.
\end{equation}

\end{lemma}

\begin{proof}
Since
\[
1-a\cdot v(p)\geq1-\gamma,
\]
the functions
\[
p\longmapsto
\frac{a}{\langle p\rangle\beta_a(p)}
=
\frac{a}{\langle p\rangle^2(1-a\cdot v(p))}
\]
and their first derivatives with respect to $p$ and $a$ are
uniformly bounded for $|a|\leq\gamma$. Since in addition $v$ is bounded and Lipschitz continuous,
the bound on $Q_{\phi,a}$ in \eqref{Bbound}
and the estimate \eqref{Bstability} follow. 
To study the diffusion coefficient $\sigma_\phi(p)=\sqrt{2\phi}\,D^{1/2}(p)$, we compute
\[
\partial_{p_k}(D^{1/2})_{ij}(p)
=
-\frac{p_k\delta_{ij}}{2\langle p\rangle^{5/2}}
+
\frac{\delta_{ik}p_j+p_i\delta_{jk}}
{\langle p\rangle^{1/2}(\langle p\rangle+1)}-
\frac{(3\langle p\rangle+1)p_ip_jp_k}
{2\langle p\rangle^{5/2}(\langle p\rangle+1)^2}.
\]
Consequently,
\[
|\partial_{p_k}(D^{1/2})_{ij}(p)|
\leq
C\langle p\rangle^{-1/2}
\leq C.
\]
Thus, $D^{1/2}$ is globally Lipschitz. Since the eigenvalues of $D(p)$ are
$\langle p\rangle$, $\langle p\rangle^{-1}$, $\langle p\rangle^{-1}$, we have
\[
\|D^{1/2}(p)\|_{\rm HS}^2
=
\mathrm{Tr}\,D(p)
=
\langle p\rangle+\frac2{\langle p\rangle}
\leq3\langle p\rangle, 
\]
and therefore the second bound in~\eqref{Bbound} follows.
Moreover,
\[
|\sqrt\phi-\sqrt\psi|
\leq
\frac{|\phi-\psi|}{2\sqrt\delta}.
\]
Writing
\[
\sigma_\phi(p)-\sigma_\psi(q)
=
\sqrt{2\phi}\,
\bigl(D^{1/2}(p)-D^{1/2}(q)\bigr)
+
\sqrt2(\sqrt\phi-\sqrt\psi)D^{1/2}(q)
\]
and using the preceding estimates gives
\eqref{sigmastability}. 
The bounds \eqref{frozenhigherderivatives}
follow by differentiating the explicit formulas
for $Q_{\phi,a}$ and $D^{1/2}$, using
$1-a\cdot v\geq1-\gamma$.
\end{proof}

\begin{theorem}[Global well-posedness of the frozen equation]
\label{frozenwellposedness}
Let $\phi$ and $a$ satisfy \eqref{frozenbounds}, and let
$G_0\in L^1(\R^3)$ satisfy
\[
G_0\geq0\ a.e.,
\qquad
\int_{\R^3}G_0\,dp=1,
\qquad
E_{G_0}<\infty.
\]
Then the Cauchy problem
\begin{equation}\label{frozenCauchy}
\begin{cases}
\displaystyle
\partial_tG
=
\phi(t)\nabla_p\cdot(D\nabla_pG)
+
\nabla_p\cdot(V_{a(t)}G),
\\[1mm]
G(0,p)=G_0(p)
\end{cases}
\end{equation}
has a weak solution
\[
G^{\phi,a}\in C([0,\infty);L^1(\R^3))
\]
such that, for every $T>0$,
\begin{equation}\label{frozenmomentbound}
\sup_{0\leq t\leq T}
E_{G^{\phi,a}(t)}
\leq
E_{G_0}+(3\Lambda_T+\mu)T.
\end{equation}
The solution is a time dependent probability density:
\begin{equation}\label{frozenmass}
G^{\phi,a}(t,\cdot)\geq0\ \text{a.e.},
\qquad
\int_{\R^3}G^{\phi,a}(t,p)\,dp=1,
\qquad t\geq0.
\end{equation}
For every $T>0$, the solution is unique among
nonnegative weak solutions in
$C([0,T];L^1(\R^3))$ with unit mass and
bounded energy on $[0,T]$.
Moreover, for every $t\geq0$, the probability measure
$G^{\phi,a}(t,p)\,dp$ is the distribution of the solution
$P_t$ of \eqref{frozenSDE} whose initial data $P_0$
has distribution $G_0(p)\,dp$.
Finally, the frozen evolution is contractive in $L^1$.
More precisely, if $G^{\phi,a}$ and $H^{\phi,a}$
correspond to initial probability densities $G_0$
and $H_0$ with finite energy, then
\begin{equation}\label{frozenL1contraction}
\|G^{\phi,a}(t)-H^{\phi,a}(t)\|_{L^1}
\leq
\|G_0-H_0\|_{L^1},
\qquad t\geq0.
\end{equation}
\end{theorem}

\begin{proof}

We divide the proof into several steps,
using standard results from the theory
of stochastic differential equations
and stochastic flows.

\medskip
\noindent
{\it Step 1. Construction by the stochastic differential equation.}

Fix $T>0$. By Lemma~\ref{frozencoefflemma} and
\eqref{frozenbounds}, the coefficients
$Q_{\phi(t),a(t)}(p)$ and $\sigma_{\phi(t)}(p)$
are globally Lipschitz in $p$, uniformly for
$t\in[0,T]$, and satisfy the linear growth
condition. Therefore, by the standard Cauchy
theory for stochastic differential
equations~\cite{Oksendal2003},
\eqref{frozenSDE} has a unique nonexplosive
strong solution for every almost surely finite
initial random variable $P_0$ independent of
the Brownian motion. Since $T>0$ is arbitrary,
the solution exists globally.
Moreover, the  regularity of the
coefficients established in
Lemma~\ref{frozencoefflemma} allows us to
apply the theory of differentiable stochastic
flows~\cite[Theorem~4.6.5]{Kunita1990}.
Consequently, the solutions of~\eqref{frozenSDE}
generate a stochastic flow $X_{s,t}$ of
$C^1$ diffeomorphisms of $\R^3$.
More precisely, let $P_t^{s,p}$ denote the
solution of~\eqref{frozenSDE} at time $t$ with deterministic initial
condition $P_s^{s,p}=p$. Then, for each
$s\leq t$ and almost every sample path of
the Brownian motion, the map
\[
X_{s,t}:\R^3\to\R^3,
\qquad X_{s,t}(p)=P_t^{s,p},
\]
is a $C^1$ diffeomorphism of $\R^3$.

Choose $P_0$ with probability distribution
$\nu_0=G_0(p)\,dp$, independently of the Brownian
motion, and denote the distribution of $P_t$
by $\nu_t$. 
The pushforward of $G_0(p)\,dp$ under $X_{0,t}$
is the random probability measure defined by
\[
(X_{0,t})_\#(G_0\,dp)(A)
=
\int_{X_{0,t}^{-1}(A)}G_0(p)\,dp,
\qquad A\subset\R^3\ \text{Borel}.
\]
Since $P_0$ is independent of the Brownian
motion, the distribution $\nu_t$ of $P_t$
satisfies
\[
\nu_t(A)=\mathbb{E}[(X_{0,t})_\#(G_0\,dp)(A)],
\qquad A\subset\R^3\ \text{Borel}.
\]
By the change of variables formula, the pushforward
measure has density
\[
g_t(p)=
G_0(X_{0,t}^{-1}(p))
\left|\det\nabla X_{0,t}^{-1}(p)\right|,
\]
and therefore
\[
\nu_t(A)=
\E\left[\int_A g_t(p)\,dp\right]
=
\int_A\E[g_t(p)]\,dp
\]
for every Borel set $A\subset\R^3$.
Thus, $\nu_t$ has density
$G^{\phi,a}(t,p)=\E[g_t(p)]$
for every $t>0$.
We set $G^{\phi,a}(0,\cdot)=G_0$.
By construction,
\[
G^{\phi,a}(t,p)\geq0\ a.e.,
\qquad
\int_{\R^3}G^{\phi,a}(t,p)\,dp=1,
\qquad t\geq0.
\]
Finally, applying It\^o's formula to
$\psi(t,P_t)$, with
$\psi\in C_c^\infty([0,T)\times\R^3)$,
and taking expectations shows, as established
in~\eqref{temporale}, that
$G^{\phi,a}$ satisfies the weak
formulation~\eqref{frozenWeakForm}.
The strong continuity of $G^{\phi,a}$
in $L^1(\R^3)$ will be established in Step 4.

\medskip
\noindent

{\it Step 2. Energy estimate.}

Let
\[
h(p)=\langle p\rangle
\]
and recall the definition~\eqref{frozenGenerator}
of the backward generator $L_t$.
A direct calculation gives
\begin{equation}\label{frozenLyapunov}
L_th
=
3\phi(t)-\mu+
\frac{\mu}{\langle p\rangle\beta_{a(t)}(p)}.
\end{equation}
Moreover,
\[
\frac1{\langle p\rangle\beta_a(p)}
=
\frac{1-|v|^2}{1-a\cdot v}
\leq
\frac{1-|v|^2}{1-|v|}
=
1+|v|
\leq2.
\]
Consequently, for $0\leq t\leq T$,
\[
-\mu\leq L_th\leq3\Lambda_T+\mu.
\]

For $R>0$, define the stopping time
\[
\tau_R
=
\inf\{t\geq0:\ |P_t|\geq R\},
\]
with the convention $\inf\varnothing=\infty$.
Since the solution of \eqref{frozenSDE}
is global and nonexplosive, we have
$\tau_R\to\infty$ almost surely as $R\to\infty$.
Applying It\^o's formula to
$h(P_{t\wedge\tau_R})$ and taking expectations,
we obtain
\[
\E h(P_{t\wedge\tau_R})
=
\E h(P_0)
+
\E\int_0^t
\mathbf1_{\{s<\tau_R\}}
L_sh(P_s)\,ds.
\]
The stopped stochastic integral has zero expectation
because its integrand is bounded on the stopped region,
and $E_{G_0}<\infty$.
Using the upper bound on $L_sh$, we deduce
\[
\E h(P_{t\wedge\tau_R})
\leq
E_{G_0}+(3\Lambda_T+\mu)t.
\]
Letting $R\to\infty$ and using Fatou's lemma gives
\[
\E h(P_t)
\leq
E_{G_0}+(3\Lambda_T+\mu)t,
\qquad 0\leq t\leq T.
\]
Since the distribution of $P_t$ has density
$G^{\phi,a}(t,\cdot)$, we conclude that
\[
E_{G^{\phi,a}(t)}
\leq
E_{G_0}+(3\Lambda_T+\mu)t,
\qquad 0\leq t\leq T.
\]
This proves \eqref{frozenmomentbound}.

\medskip
\noindent
{\it Step 3. The $L^1$ contraction property.}

For $0\leq s\leq t$, let $P_t^{s,q}$ denote the
solution of~\eqref{frozenSDE} at time $t$ with
deterministic initial condition $P_s^{s,q}=q$.
Define the transition probability
\[
K_{s,t}(q,A)
=
\mathbb P(P_t^{s,q}\in A),
\]
for every Borel set $A\subset\R^3$. 
The stochastic flow constructed in Step 1
shows that, for every $f\in L^1(\R^3)$,
the measure
\[
A\longmapsto
\int_{\R^3}K_{s,t}(q,A)f(q)\,dq
\]
is absolutely continuous with respect to
Lebesgue measure. Thus, we may define the
forward propagator $U(t,s)$ by
\begin{equation}\label{defU}
\int_A[U(t,s)f](p)\,dp
=
\int_{\R^3}K_{s,t}(q,A)f(q)\,dq.
\end{equation}
By the Markov property,
\[
U(t,r)U(r,s)=U(t,s),
\qquad 0\leq s\leq r\leq t.
\]
In particular,
\[
G^{\phi,a}(t,\cdot)
=
U(t,s)G^{\phi,a}(s,\cdot).
\]
By positivity and linearity, 
$
|U(t,s)f|\leq U(t,s)|f|.
$
Moreover, conservation of total probability gives
\[
\int_{\R^3}U(t,s)|f|\,dp
=
\int_{\R^3}|f|\,dp.
\]
Hence 
\[
\|U(t,s)f\|_{L^1}
\leq
\|f\|_{L^1}.
\]
Applying this inequality to $f=G_0-H_0$
proves~\eqref{frozenL1contraction}.

\medskip
\noindent
{\it Step 4. Strong continuity in $L^1$.}

We prove that, for every $f\in L^1(\R^3)$,
the map $t\mapsto U(t,0)f$ is continuous
from $[0,\infty)$ into $L^1(\R^3)$.
Let $X_{s,t}$ be the stochastic flow of
$C^1$ diffeomorphisms constructed in Step 1.
For $f\in L^1(\R^3)$, define
\begin{equation}\label{push}
f_{0,t}(p)
=
f(X_{0,t}^{-1}(p))
\left|\det\nabla X_{0,t}^{-1}(p)\right|.
\end{equation}
By the change of variables formula,
$f_{0,t}(p)\,dp$ is the pushforward of
$f(p)\,dp$ under $X_{0,t}$. In particular,
\begin{equation}\label{intf}
\|f_{0,t}\|_{L^1}=\|f\|_{L^1}
\end{equation}
almost surely. Moreover, the definition
of the transition probabilities and
Fubini's theorem give
\begin{equation}\label{Urep}
U(t,0)f=\E[f_{0,t}]
\end{equation}
as an equality in $L^1(\R^3)$. 

Suppose first that $f\in C_c^\infty(\R^3)$.
By the regularity theory of differentiable
stochastic flows
\cite[Theorem~4.6.5]{Kunita1990},
we may choose a version of the flow such
that, almost surely, $X_{0,t}$, its inverse,
and their Jacobian matrices are continuous
in $(t,p)$.
Consequently, for almost every sample path of the Brownian motion,
\[
f_{0,t}(p)\rightarrow f_{0,s}(p)
\]
locally uniformly in $p$ as $t\to s$. 
Let $K=\operatorname{supp}f$.
Since $(t,p)\mapsto X_{0,t}(p)$ is
continuous and $K$ is compact, for almost
every path of the Brownian motion, there exist $\varepsilon>0$
and $R>0$ such that $X_{0,t}(K)\subset B_R$ 
whenever $t\geq0$ and $|t-s|<\varepsilon$.
Thus, $f_{0,t}$ and $f_{0,s}$ are supported
in the same bounded set $B_R$ for all such $t$.
It follows that
\[
\|f_{0,t}-f_{0,s}\|_{L^1}
=
\int_{B_R}|f_{0,t}(p)-f_{0,s}(p)|\,dp
\leq
|B_R|\sup_{p\in B_R}
|f_{0,t}(p)-f_{0,s}(p)|
\rightarrow0
\]
almost surely as $t\to s$. 
Since \eqref{intf} implies $\|f_{0,t}-f_{0,s}\|_{L^1}
\leq2\|f\|_{L^1}$, 
the dominated convergence theorem gives
\[
\E\|f_{0,t}-f_{0,s}\|_{L^1}\to0.
\]
Using \eqref{Urep}, we obtain
\[
\|U(t,0)f-U(s,0)f\|_{L^1}
\leq
\E\|f_{0,t}-f_{0,s}\|_{L^1}
\rightarrow0.
\]
Finally, the result extends to every
$f\in L^1(\R^3)$ by the density of
$C_c^\infty(\R^3)$ in $L^1(\R^3)$
and the contraction property established
in Step 3.
In particular, taking $f=G_0$ proves
\[
G^{\phi,a}\in C([0,\infty);L^1(\R^3)),
\]
as required.

\medskip
\noindent
{\it Step 5. Uniqueness.}

The solution $G=G^{\phi,a}$ constructed in Step 1
is the probability density of the solution
$P_t$ of~\eqref{frozenSDE}, whose initial
distribution has density $G_0$.
In particular, $P_t$ solves the martingale
problem associated with the
generator~\eqref{GenSDE}.
We recall that a continuous stochastic process
$P_t$ solves the martingale problem associated
with $L_t$ if, for every
$f\in C_c^\infty(\R^3)$, the process
\[
M_t^f
=
f(P_t)-f(P_0)
-\int_0^t L_sf(P_s)\,ds
\]
is a martingale with respect to the natural
filtration of $P_t$.
To prove uniqueness in
Theorem~\ref{frozenwellposedness},
we use the superposition principle
of~\cite{Trevisan2016} to show that any given weak solution $\widetilde G$ of~\eqref{frozenCauchy}
in the class stated in the theorem gives rise to a solution of the martingale problem.

Assume $\widetilde G(0,\cdot )=G_0$. Fix $T>0$ and set
\[
\widetilde\nu_t(dp)=\widetilde G(t,p)\,dp.
\]
Recall that the frozen equation can be written
in the Fokker--Planck--Kolmogorov form
\[
\partial_t\widetilde G
=\sum_{i,j=1}^3
\partial_{p_i}\partial_{p_j}
(A_{ij}\widetilde G)
-
\nabla_p\cdot(Q\widetilde G),
\]
where
\[
A(t,p)=\phi(t)D(p),
\qquad
Q(t,p)=Q_{\phi(t),a(t)}(p).
\]
By Lemma~\ref{frozencoefflemma},
\[
\|A(t,p)\|\leq C_T\langle p\rangle,
\qquad
|Q(t,p)|\leq C_T.
\]
Since $\widetilde G$ has unit mass and
locally bounded energy,
\[
\int_0^T\int_{\R^3}
\bigl(\|A(t,p)\|+|Q(t,p)|\bigr)
\,d\widetilde\nu_t(p)\,dt<\infty.
\]
Moreover, the $L^1$ continuity of
$\widetilde G$ implies that
$t\mapsto\widetilde\nu_t$ is narrowly
continuous.
We may therefore apply the superposition
principle~\cite[Theorem~2.5]{Trevisan2016}.
It follows that there exists a continuous
stochastic process $\widetilde P_t$
which solves the martingale problem
associated with the generator
\[
L_tf
=
Q\cdot\nabla_pf+A:\nabla_p^2f,
\]
and whose time marginal distributions are $\widetilde\nu_t$, $0\leq t\leq T$. 
Since the coefficients of~\eqref{frozenSDE}
are globally Lipschitz in $p$ and satisfy
the linear growth condition, the associated
martingale problem is well posed.
Since $\widetilde P_0$ and $P_0$
have the same distribution $G_0(p)\,dp$,
the processes $\widetilde P_t$ and $P_t$
have also the same distribution.
Consequently,
\[
\widetilde\nu_t(dp)
=
G(t,p)\,dp,
\qquad 0\leq t\leq T.
\]
Since $T>0$ is arbitrary, uniqueness follows.
\end{proof}

The following Wasserstein estimates will be used
in the proof of local existence for the nonlinear
equation. The first estimate controls the time
continuity of the nonlinear functionals
$(\Phi_G,j_G)$, while the second provides the
stability with respect to the frozen coefficients
needed to establish contraction of the fixed-point
map.

\begin{lemma}[Frozen Wasserstein estimates]
\label{frozenestimates}
Let $T>0$ and assume that
$\phi,\psi\in C([0,T])$ and
$a,b\in C([0,T];\R^3)$ satisfy
\[
0<\delta\leq\phi,\psi\leq\Lambda,
\qquad
|a|,|b|\leq\gamma<1.
\]
Let $G_0\in L^1(\R^3)$ be a nonnegative
probability density with $E_{G_0}<\infty$,
and denote by $G^{\phi,a}$ and $G^{\psi,b}$
the corresponding solutions of the frozen
Fokker--Planck equation~\eqref{frozenFP}.
There exists a constant $C_T$, depending only
on $T,\delta,\Lambda,\gamma,\mu$ and
$E_{G_0}$, such that
\begin{equation}\label{frozenW1time}
W_1\bigl(
G^{\phi,a}(t),G^{\phi,a}(s)
\bigr)
\leq C_T|t-s|^{1/2},
\qquad 0\leq s\leq t\leq T,
\end{equation}
and
\begin{equation}\label{frozenW1stability}
W_1\bigl(
G^{\phi,a}(t),G^{\psi,b}(t)
\bigr)
\leq
C_T\sqrt t\,
\|(\phi,a)-(\psi,b)\|_\infty,
\qquad 0\leq t\leq T.
\end{equation}
\end{lemma}

\begin{proof}
We use the stochastic representation established
in Theorem~\ref{frozenwellposedness}.
Throughout the proof we use the coefficient
estimates of Lemma~\ref{frozencoefflemma},
with constants depending only on
$\delta,\Lambda,\gamma$ and $\mu$.

Let $P_t$ be the solution of~\eqref{frozenSDE} corresponding to
$(\phi,a)$, with initial data $P_0$ having probability
distribution $\nu_0=G_0(p)\,dp$.
The energy estimate in
Theorem~\ref{frozenwellposedness} gives
\[
\sup_{0\leq t\leq T}
\E\langle P_t\rangle
\leq K_T,
\]
where
\[
K_T=E_{G_0}+(3\Lambda+\mu)T.
\]
For $0\leq s\leq t\leq T$, the SDE~\eqref{frozenSDE} gives
\[
P_t-P_s
=
\int_s^t Q_{\phi(r),a(r)}(P_r)\,dr
+
\int_s^t\sigma_{\phi(r)}(P_r)\,dB_r.
\]
By the Cauchy--Schwarz inequality and
It\^o's isometry,
\[
\E|P_t-P_s|^2
\leq
2(t-s)\int_s^t
\E|Q_{\phi(r),a(r)}(P_r)|^2\,dr+
2\int_s^t
\E\|\sigma_{\phi(r)}(P_r)\|_{\mathrm{HS}}^2\,dr\leq
C_T|t-s|.
\]
This estimate does not require
$\E|P_0|^2<\infty$, since it concerns
only the increment $P_t-P_s$.
Since the joint distribution of $(P_t,P_s)$ is a coupling
of $G^{\phi,a}(t)$ and $G^{\phi,a}(s)$,
\[
W_1\bigl(
G^{\phi,a}(t),G^{\phi,a}(s)
\bigr)
\leq\E|P_t-P_s|\leq
\bigl(\E|P_t-P_s|^2\bigr)^{1/2}\leq C_T|t-s|^{1/2}.
\]
This proves \eqref{frozenW1time}.
We next prove \eqref{frozenW1stability}.
Let $P_t$ and $\widetilde P_t$ be the solutions
corresponding to $(\phi,a)$ and $(\psi,b)$,
respectively. We construct them with the same
initial random variable, having distribution $G_0(p)\,dp$,
and driven by the same Brownian motion.
Set
\[
Z_t=P_t-\widetilde P_t,
\qquad Z_0=0,
\]
and
\[
\eta=
\|(\phi,a)-(\psi,b)\|_\infty.
\]
For $R>0$, define the stopping time
\[
\tau_R
=
\inf\left\{
t\geq0:
|P_t|+|\widetilde P_t|\geq R
\right\}.
\]
Applying It\^o's formula to
$|Z_{t\wedge\tau_R}|^2$, taking expectations,
and using the coefficient estimates above
together with Young's inequality, we obtain
\[
\E|Z_{t\wedge\tau_R}|^2
\leq
C\int_0^t
\E|Z_{s\wedge\tau_R}|^2\,ds+
C\eta^2\int_0^t
\left(
1+\E\langle\widetilde P_s\rangle
\right)ds.
\]
The energy estimate gives
\[
\sup_{0\leq s\leq T}
\E\langle\widetilde P_s\rangle
\leq K_T.
\]
Consequently,
\[
\E|Z_{t\wedge\tau_R}|^2
\leq
C\int_0^t
\E|Z_{s\wedge\tau_R}|^2\,ds
+
C_Tt\eta^2.
\]
By Gronwall's inequality,
\[
\E|Z_{t\wedge\tau_R}|^2
\leq C_Tt\eta^2.
\]
Since both SDE solutions are global and
nonexplosive, $\tau_R\to\infty$ almost surely
as $R\to\infty$. Fatou's lemma therefore yields
\[
\E|Z_t|^2
\leq C_Tt\eta^2.
\]
Finally, the joint distribution of
$(P_t,\widetilde P_t)$ is a coupling of the
two frozen solutions. Hence
\[
W_1\bigl(
G^{\phi,a}(t),G^{\psi,b}(t)
\bigr)
\leq
\E|P_t-\widetilde P_t|
\leq
\bigl(\E|Z_t|^2\bigr)^{1/2}\leq C_T\sqrt t\,\eta.
\]
This proves \eqref{frozenW1stability}.
\end{proof}

%===========================================================
\section{Proof of Theorem~\ref{mainthm}}
%===========================================================

We are now ready to prove the first main theorem of this
paper. We split the proof into several subsections.

%===========================================================
\subsection{Local existence}
%===========================================================

We first use the frozen problem studied in
Section~\ref{frozensec} to construct a local solution of the
nonlinear equation. Recall that
\[
j_G=\int_{\R^3}v(p)G(p)\,dp
\]
and
\[
\Phi_G
=
\frac{\mu}{3}
\left[
1-
\int_{\R^3}
\frac{G(p)}
{\langle p\rangle\beta_G(p)}\,dp
\right],
\qquad
\beta_G(p)
=
\langle p\rangle
\bigl(1-j_G\cdot v(p)\bigr).
\]
Let $G_0$ satisfy
\[
G_0\in L^1(\R^3),\qquad
G_0\geq0,\qquad
\int_{\R^3}G_0\,dp=1,\qquad
E_{G_0}<\infty.
\]
Set
\[
a_0=j_{G_0},
\qquad
\phi_0=\Phi_{G_0}.
\]
By Lemmas~\ref{estimatesBeta}
and~\ref{Phiboundslemma},
\[
|a_0|<1,
\qquad
0<\phi_0\leq\frac{\mu}{3}.
\]
We define
\begin{equation}\label{localconstants}
\delta=\frac{\phi_0}{2},
\qquad
\gamma=\frac{1+|a_0|}{2},
\qquad
\Lambda=\frac{\mu}{3}.
\end{equation}
Thus
\[
0<\delta<\phi_0\leq\Lambda,
\qquad
|a_0|<\gamma<1.
\]
For $0<T\leq1$, let
\begin{align}
\mathcal K_T
=
\bigl\{(\phi,a)\in
C([0,T];\R\times\R^3):\;&
\delta\leq\phi(t)\leq\Lambda,\quad
|a(t)|\leq\gamma,
\nonumber\\
&
\phi(0)=\phi_0,\quad a(0)=a_0
\bigr\}.
\label{fixedpointset}
\end{align}
We equip $\mathcal K_T$ with the norm
\[
\|(\phi,a)\|_\infty
=
\sup_{0\leq t\leq T}
\bigl(|\phi(t)|+|a(t)|\bigr).
\]
Then $\mathcal K_T$ is a nonempty closed subset of the Banach
space $C([0,T];\R\times\R^3)$
and is therefore complete. 

For every $(\phi,a)\in\mathcal K_T$, let
$G^{\phi,a}$ denote the unique solution
on $[0,T]$ of the frozen equation
\begin{equation}\label{localfrozen}
\partial_tG^{\phi,a}
=
\phi(t)\nabla_p\cdot
\bigl(D\nabla_pG^{\phi,a}\bigr)
+
\nabla_p\cdot
\bigl(V_aG^{\phi,a}\bigr),
\qquad
G^{\phi,a}(0)=G_0.
\end{equation}
By Theorem~\ref{frozenwellposedness},
\[
G^{\phi,a}\in C([0,T];L^1(\R^3)),
\qquad
G^{\phi,a}\geq0,
\qquad
\int_{\R^3}G^{\phi,a}(t,p)\,dp=1.
\]
We define
\begin{equation}\label{fixedpointmap}
\mathcal F(\phi,a)(t)
=
\left(
\Phi_{G^{\phi,a}(t)},
j_{G^{\phi,a}(t)}
\right).
\end{equation}
We first record estimates whose constants can be chosen independently
of $T\leq1$. Indeed, the constants in
Theorem~\ref{frozenwellposedness}
and Lemma~\ref{frozenestimates} depend only on the upper
bound for the time interval and on the parameters $
\delta,\Lambda,\gamma,\mu,E_{G_0}$. 
Extend $\phi$ and $a$ constantly from $T$ to $1$.
The extended coefficients satisfy the same bounds
$\delta\leq\phi\leq\Lambda$ and $|a|\leq\gamma$.
Hence, applying Theorem~\ref{frozenwellposedness}
and Lemma~\ref{frozenestimates}
with upper time horizon $1$, there are constants
$K$ and $C_0$, independent of $T\leq1$ and of
$(\phi,a)\in\mathcal K_T$, such that
\begin{equation}\label{localmomentuniform}
\sup_{0\leq t\leq T}
E_{G^{\phi,a}(t)}
\leq K,
\end{equation}
where one may take
\[
K=E_{G_0}+2\mu,
\]
and
\begin{equation}\label{localtimeW1}
W_1\bigl(
G^{\phi,a}(t),G^{\phi,a}(s)
\bigr)
\leq
C_0|t-s|^{1/2}.
\end{equation}
Moreover, for
$(\phi,a),(\psi,b)\in\mathcal K_T$,
\begin{equation}\label{localfrozenstability}
W_1\bigl(
G^{\phi,a}(t),G^{\psi,b}(t)
\bigr)
\leq
C_0\sqrt t\,
\|(\phi,a)-(\psi,b)\|_\infty.
\end{equation}
By Lemma~\ref{estimatesBeta} and
\eqref{localmomentuniform},
\[
|j_{G^{\phi,a}(t)}|
\leq
\sqrt{
1-\frac{1}{E_{G^{\phi,a}(t)}^2}
}
\leq
\sqrt{1-\frac1{K^2}}
=:\gamma_*<1.
\]
Consequently, all the probability densities occurring in the
fixed point argument satisfy a common strict bound
\[
|j_G|\leq\gamma_*<1.
\]
By Lemma~\ref{momentLipschitz}, there exists a constant $C_1$,
depending only on $\gamma_*$ and $\mu$, such that
\begin{equation}\label{functionalW1}
|\Phi_G-\Phi_H|+|j_G-j_H|
\leq
C_1W_1(G,H)
\end{equation}
for all probability densities $G,H$ arising in the construction.

We can now prove local existence.

\begin{theorem}[Local existence]
\label{localexistencetheorem}
Let
\[
G_0\in L^1(\R^3),\qquad
G_0\geq0\ a.e.,\qquad
\int_{\R^3}G_0\,dp=1,\qquad
E_{G_0}<\infty.
\]
Then there exists $T_0>0$ such that the
Cauchy problem for~\eqref{mainPDE} with initial data $G_0$ has a weak solution on $[0,T_0]$.
Moreover,
\[
G\in C([0,T_0];L^1(\R^3)),
\qquad
G(t,\cdot)\geq0\ a.e.,
\qquad
\int_{\R^3}G(t,p)\,dp=1,
\]
and
\[
\sup_{0\leq t\leq T_0}E_{G(t)}<\infty.
\]
In addition, if, for some constants
\[
0<\delta_0\leq\frac{\mu}{3},
\qquad
0\leq\gamma_0<1,
\qquad
1\leq K_0<\infty,
\]
there holds
\[
\Phi_{G_0}\geq\delta_0,
\qquad
|j_{G_0}|\leq\gamma_0,
\qquad
E_{G_0}\leq K_0,
\]
then $T_0$ may be chosen depending only on
$\delta_0,\gamma_0,K_0$ and $\mu$.
\end{theorem}

\begin{proof}
We claim that, for $T>0$ sufficiently small, the map
$\mathcal F$ defined by \eqref{fixedpointmap} is a contraction of
$\mathcal K_T$ into itself.
We first prove that $\mathcal F$ maps $\mathcal K_T$ into itself.
Since $G^{\phi,a}(0)=G_0$, we have
\[
\mathcal F(\phi,a)(0)
=
(\phi_0,a_0).
\]
Moreover, by \eqref{functionalW1} and
\eqref{localtimeW1},
\begin{equation}
\left|
\Phi_{G^{\phi,a}(t)}-\phi_0
\right|
+
\left|
j_{G^{\phi,a}(t)}-a_0
\right|\leq
C_1W_1\bigl(
G^{\phi,a}(t),G_0
\bigr)
\leq
C_1C_0\sqrt t.
\label{selfmapestimate}
\end{equation}
In particular,
\[
\Phi_{G^{\phi,a}(t)}
\geq
\phi_0-C_1C_0\sqrt T.
\]
Choose $T$ so small that
\begin{equation}\label{selfmapT}
C_1C_0\sqrt T
\leq
\min\left\{
\frac{\phi_0}{2},
\gamma-|a_0|
\right\}.
\end{equation}
Then
\[
\Phi_{G^{\phi,a}(t)}
\geq
\frac{\phi_0}{2}
=
\delta
\]
and
\[
|j_{G^{\phi,a}(t)}|
\leq
|a_0|+C_1C_0\sqrt T
\leq\gamma.
\]
On the other hand, Lemma~\ref{Phiboundslemma} gives
\[
\Phi_{G^{\phi,a}(t)}
\leq
\frac{\mu}{3}
=
\Lambda.
\]
Thus, the pointwise bounds defining $\mathcal K_T$ are satisfied.
Moreover, again by \eqref{functionalW1} and
\eqref{localtimeW1},
\[
\left|
\Phi_{G^{\phi,a}(t)}
-
\Phi_{G^{\phi,a}(s)}
\right|
+
\left|
j_{G^{\phi,a}(t)}
-
j_{G^{\phi,a}(s)}
\right|\leq
C_1C_0|t-s|^{1/2}.
\]
Hence, the functions defining the image of $\mathcal F$
are continuous and therefore
$\mathcal F(\phi,a)\in\mathcal K_T$.
We next prove that $\mathcal F$ is a contraction. Let
\[
(\phi,a),(\psi,b)\in\mathcal K_T.
\]
By \eqref{functionalW1} and
\eqref{localfrozenstability},
\begin{align*}
\left|
\Phi_{G^{\phi,a}(t)}
-
\Phi_{G^{\psi,b}(t)}
\right|
+
\left|
j_{G^{\phi,a}(t)}
-
j_{G^{\psi,b}(t)}
\right|
&\leq
C_1
W_1\bigl(
G^{\phi,a}(t),G^{\psi,b}(t)
\bigr)\\
&\leq
C_1C_0\sqrt t\,
\|(\phi,a)-(\psi,b)\|_\infty.
\end{align*}
Taking the supremum over $0\leq t\leq T$ gives
\[
\|
\mathcal F(\phi,a)-\mathcal F(\psi,b)
\|_\infty
\leq
C_1C_0\sqrt T\,
\|(\phi,a)-(\psi,b)\|_\infty.
\]
Shrinking $T$ further, if necessary, so that
\begin{equation}\label{contractionT}
C_1C_0\sqrt T<1,
\end{equation}
the map $\mathcal F$ is a strict contraction on the complete
metric space $\mathcal K_T$.
Banach's fixed-point theorem therefore gives a unique
$(\phi,a)\in\mathcal K_T$ such that
\[
\mathcal F(\phi,a)=(\phi,a).
\]
If $G=G^{\phi,a}$, then, by the definition of $\mathcal F$,
\[
\phi(t)=\Phi_{G(t)},
\qquad
a(t)=j_{G(t)}.
\]
Consequently \eqref{localfrozen} becomes
\[
\partial_tG
=
\Phi_G\nabla_p\cdot(D\nabla_pG)
+
\nabla_p\cdot(V_{j_G}G),
\]
which is precisely the nonlinear equation~\eqref{mainPDE}.
The properties
\[
G\in C([0,T];L^1(\R^3)),
\qquad
G(t,\cdot)\geq0\ a.e.,
\qquad
\int_{\R^3}G(t,p)\,dp=1
\]
follow from Theorem~\ref{frozenwellposedness}, while
\[
\sup_{0\leq t\leq T}E_{G(t)}\leq K
\]
follows from \eqref{localmomentuniform}.
This proves local existence.
Finally, suppose that
\[
\phi_0=\Phi_{G_0}\geq\delta_0,
\qquad
|a_0|=|j_{G_0}|\leq\gamma_0<1,
\qquad
E_{G_0}\leq K_0.
\]
Then, by \eqref{localconstants},
\[
\delta=\frac{\phi_0}{2}\geq\frac{\delta_0}{2},
\qquad
\gamma-|a_0|
=
\frac{1-|a_0|}{2}
\geq
\frac{1-\gamma_0}{2}.
\]
Moreover, the constants $K$ and $C_0$ above
can be chosen to depend only on
$\delta_0,\gamma_0,K_0$ and $\mu$.
In particular, we may take
\[
K=K_0+2\mu,
\]
so that
\[
\gamma_*=\sqrt{1-\frac1{K^2}}<1
\]
and the constant $C_1$ can also be chosen uniformly
in terms of $\delta_0,\gamma_0,K_0$ and $\mu$.
Therefore, the smallness conditions
\eqref{selfmapT} and \eqref{contractionT}
are satisfied for a common positive time $T_0$
depending only on $\delta_0,\gamma_0,K_0$ and $\mu$.
\end{proof}

Global existence and uniqueness will be
established below.

%===========================================================
\subsection{Conservation laws}
%===========================================================

We next prove conservation of energy and momentum.
For smooth solutions of~\eqref{mainPDE}, this property
was established in~\cite{Calogero}.
The argument below extends this result to weak
solutions with finite initial energy.

Let $G$ be a weak solution on $[0,T]$ with initial data $G_0$
satisfying
\[
\int_{\R^3}\langle p\rangle G_0(p)\,dp<\infty.
\]
For smooth functions $f(p)$, introduce the nonlinear operator
\[
L_Gf
=
\Phi_G\nabla_p\cdot(D\nabla_pf)
-
W_G\cdot\nabla_pf.
\]
Formally, the weak equation gives
\[
\frac{d}{dt}
\int_{\R^3}f(p)G(t,p)\,dp
=
\int_{\R^3}L_Gf(p)G(t,p)\,dp.
\]
We justify this identity for the unbounded observables
$\langle p\rangle$ and $p_i$.

%-----------------------------------------------------------
\subsubsection*{Conservation of energy}
%-----------------------------------------------------------

Set
\[
e(p)=\langle p\rangle.
\]
A direct calculation gives
\begin{equation}\label{generatorenergy}
L_Ge
=
3\Phi_G-\mu
+
\frac{\mu}{\langle p\rangle\beta_G}.
\end{equation}
By~\eqref{betanoen} and~\eqref{Phinoen},
\begin{equation}\label{generatorenergybound}
|L_Ge(t,p)|\leq2\mu,
\qquad
0\leq t\leq T,\quad p\in\R^3.
\end{equation}
To justify the use of $e$ as a test function,
we first establish a bound on the energy of $G$.

Choose a smooth, nondecreasing, concave
function $\eta:[0,\infty)\to[0,\infty)$ satisfying
\[
\eta(r)=r\quad(0\leq r\leq1),
\qquad
\eta'(r)=0\quad(r\geq2),
\qquad
0\leq\eta'\leq1.
\]
For $R\geq1$, define
\[
e_R(p)
=
R\eta\left(\frac{e(p)}R\right).
\]
Then
\[
0\leq e_R\leq e,
\qquad
e_R\rightarrow e
\]
pointwise as $R\to\infty$. We compute
\[
L_Ge_R
=
\eta'\left(\frac eR\right)L_Ge
+
\frac{\Phi_G}{R}
\eta''\left(\frac eR\right)
\nabla_pe\cdot D\nabla_pe.
\]
Since $\eta''\leq0$, $\Phi_G\geq0$, and
\[
\nabla_pe\cdot D\nabla_pe
=
\frac{|p|^2}{\langle p\rangle}\geq0,
\]
we obtain
\[
L_Ge_R\leq2\mu.
\]
Moreover, $L_Ge_R$ is bounded uniformly
in $R$. Indeed, $\eta''(e/R)$ vanishes
outside the region $R\leq e\leq2R$,
where
\[
\frac1R\nabla_pe\cdot D\nabla_pe
\leq2.
\]
Although $e_R$ is bounded, it is not
compactly supported. Let
$\chi\in C_c^\infty([0,\infty))$ satisfy
\[
0\leq\chi\leq1,\qquad
\chi(r)=1\quad(r\leq1),\qquad
\chi(r)=0\quad(r\geq2),
\]
and set
\[
\chi_N(p)=\chi\left(\frac{|p|}{N}\right).
\]
For fixed $R$ and $N$ sufficiently large,
the function $e_R$ is constant on the
support of $\nabla_p\chi_N$.
Expanding the generator gives
\[
L_G(\chi_Ne_R)
=
\chi_NL_Ge_R
+
\Phi_Ge_R\nabla_p\cdot(D\nabla_p\chi_N)
-
e_RW_G\cdot\nabla_p\chi_N,
\]
where the cross term vanishes.
Since $\chi_N$ is radial,
\[
D\nabla_p\chi_N
=
\langle p\rangle\nabla_p\chi_N,
\]
and
\[
|\nabla_p\chi_N|\leq\frac{C}{N},
\qquad
\left|
\nabla_p\cdot(D\nabla_p\chi_N)
\right|
\leq\frac{C}{N}.
\]
Moreover, \eqref{betanoen} implies
\[
|W_G(t,p)|\leq3\mu,
\]
while $0\leq e_R\leq CR$.
Consequently,
\[
\left|
L_G(\chi_Ne_R)-\chi_NL_Ge_R
\right|
\leq C\mu\frac{R}{N}.
\]
Using $\chi_Ne_R$ in the integrated weak
formulation and letting $N\to\infty$,
we obtain
\[
\int_{\R^3}e_R(p)G(t,p)\,dp
=
\int_{\R^3}e_R(p)G_0(p)\,dp+
\int_0^t\int_{\R^3}
L_Ge_R(p)G(s,p)\,dp\,ds.
\]
Since $L_Ge_R\leq2\mu$ and $G$ has unit
mass,
\[
\int_{\R^3}e_R(p)G(t,p)\,dp
\leq
\int_{\R^3}e_R(p)G_0(p)\,dp+2\mu t.
\]
Letting $R\to\infty$, Fatou's lemma
and dominated convergence yield
\begin{equation}\label{firstmomentbound}
\int_{\R^3}\langle p\rangle G(t,p)\,dp
\leq
\int_{\R^3}\langle p\rangle G_0(p)\,dp
+2\mu t,
\qquad 0\leq t\leq T.
\end{equation}
Thus, the energy remains finite in any compact interval of time. 
We can now establish energy conservation.
For $\chi_R(p)=\chi(|p|/R)$, we have
\begin{equation}\label{energycutoffidentity}
L_G(\chi_Re)
=
\chi_RL_Ge
+
\Phi_Ge\,\nabla_p\cdot(D\nabla_p\chi_R)
+
2\Phi_G\nabla_p\chi_R\cdot D\nabla_pe
-
eW_G\cdot\nabla_p\chi_R.
\end{equation}
Since $D\nabla_pe=p$, the cutoff terms
are supported in
\[
A_R=\{R\leq|p|\leq2R\}.
\]
On this annulus,
\[
|\nabla_p\chi_R|\leq\frac{C}{R},
\qquad
\left|
\nabla_p\cdot(D\nabla_p\chi_R)
\right|
\leq\frac{C}{R},
\]
and $\langle p\rangle\leq CR$.
Using $\Phi_G\leq\mu/3$ and
$|W_G|\leq3\mu$, we obtain
\[
\left|
L_G(\chi_Re)-\chi_RL_Ge
\right|
\leq
C\mu\mathbf1_{A_R}(p),
\]
where $C$ is independent of $R$ and $T$.
These errors vanish after integration
against $G$ over $[0,t]\times\R^3$
as $R\to\infty$, by dominated convergence.
Moreover, \eqref{firstmomentbound}
ensures that the energy is finite,
while \eqref{generatorenergybound}
provides an integrable bound for
the generator term.
Passing to the limit in the integrated
weak formulation yields
\[
\int_{\R^3}e(p)G(t,p)\,dp
=
\int_{\R^3}e(p)G_0(p)\,dp+
\int_0^t
\left[
3\Phi_G(s)-\mu
+
\mu\int_{\R^3}
\frac{G(s,p)}
{\langle p\rangle\beta_G(s,p)}\,dp
\right]ds.
\]
By the definition of $\Phi_G$, the
integrand in square brackets vanishes
identically. Hence
\begin{equation}\label{energyconservation}
\int_{\R^3}\langle p\rangle G(t,p)\,dp
=
\int_{\R^3}\langle p\rangle G_0(p)\,dp
=:E,
\qquad 0\leq t\leq T.
\end{equation}

%-----------------------------------------------------------
\subsubsection*{Conservation of momentum}
%-----------------------------------------------------------

Since $|p|\leq\langle p\rangle$, energy conservation
implies that the momentum
\[
P_G(t)=\int_{\R^3}p\,G(t,p)\,dp
\]
is well defined.
For $i=1,2,3$, we have
\[
L_Gp_i
=
(3\Phi_G-\mu)v_i
+
\frac{\mu(j_G)_i}{\langle p\rangle\beta_G}.
\]
By \eqref{betanoen} and \eqref{Phinoen},
\[
|L_Gp_i|\leq3\mu.
\]
Using $\chi_Rp_i$ as a test function, the same
cutoff estimates as in the energy argument give
\[
\left|
L_G(\chi_Rp_i)-\chi_RL_Gp_i
\right|
\leq
C\mu\mathbf1_{A_R}(p),
\]
where $C$ is independent of $R$ and $T$.
The cutoff errors vanish after integration
against $G$ over $[0,t]\times\R^3$
as $R\to\infty$. Moreover, energy conservation
and the bound on $L_Gp_i$ justify passing
to the limit in the remaining terms.
We obtain
\begin{align*}
\int_{\R^3}p_iG(t,p)\,dp
={}&
\int_{\R^3}p_iG_0(p)\,dp\\
&+
\int_0^t
\left[
(3\Phi_G(s)-\mu)(j_G(s))_i
+
\mu(j_G(s))_i
\int_{\R^3}
\frac{G(s,p)}
{\langle p\rangle\beta_G(s,p)}\,dp
\right]ds.
\end{align*}
By the definition of $\Phi_G$, the integrand
in square brackets vanishes identically.
Therefore
\begin{equation}\label{momentumconservation}
\int_{\R^3}p\,G(t,p)\,dp
=
\int_{\R^3}p\,G_0(p)\,dp
=:P,
\qquad 0\leq t\leq T.
\end{equation}
Thus every weak solution with finite initial
energy conserves energy and momentum.

%===========================================================
\subsection{Propagation of energy tails}
%===========================================================

We now prove that the energy tails of weak solutions
remain uniformly integrable on every finite time interval.

For $\phi\geq0$ and $|a|<1$, let
$L_{\phi,a}$ denote the frozen generator on the
right hand side of~\eqref{frozenGenerator},
that is,
\begin{equation}\label{frozenGenerator2}
L_{\phi,a}\psi
=
\phi\nabla_p\cdot(D\nabla_p\psi)
-
V_a\cdot\nabla_p\psi.
\end{equation}
The following lemma provides an estimate
for the frozen generator that will be used
to control the energy tails of solutions
to the nonlinear equation. An important feature
of the estimate is that it is uniform in $a$,
even as $|a|$ approaches one.
\begin{lemma}\label{energytailtestlemma}
Let $\eta\in C^\infty([0,\infty))$ be
nondecreasing and satisfy
\[
0\leq\eta\leq1,
\qquad
\eta(s)=0\quad(s\leq1),
\qquad
\eta(s)=1\quad(s\geq2).
\]
For $R>1$, define
\[
H_R(p)=
\langle p\rangle
\eta\left(\frac{\langle p\rangle}{R}\right).
\]
Then
\begin{equation}\label{energytailtestbounds}
\langle p\rangle
\mathbf1_{\{\langle p\rangle\geq2R\}}
\leq H_R(p)
\leq
\langle p\rangle
\mathbf1_{\{\langle p\rangle>R\}}.
\end{equation}
Moreover, for every $\Lambda>0$, there exists
a constant $C>0$, depending only on $\eta$,
such that
\begin{equation}\label{energytailgenerator}
|L_{\phi,a}H_R(p)|
\leq
C(\Lambda+\mu)
\mathbf1_{\{\langle p\rangle>R\}},
\end{equation}
for all $0\leq\phi\leq\Lambda$ and $|a|<1$.
The constant $C$ is independent of
$R,\phi,a$ and $p$.
\end{lemma}

\begin{proof}
Set $h=\langle p\rangle$. A direct computation gives
\[
L_{\phi,a}h
=
3\phi-\mu+
\frac{\mu}{h\beta_a}.
\]
Moreover,
\[
\frac1{h\beta_a}
=
\frac{1-|v|^2}{1-a\cdot v}
\leq
\frac{1-|v|^2}{1-|v|}
=
1+|v|
\leq2.
\]
Consequently,
\[
|L_{\phi,a}h|\leq3\Lambda+\mu.
\]
Write $
H_R(p)=g_R(h)$, $g_R(h)=h\eta(h/R)$. 
Then
\[
g_R'(h)
=
\eta(h/R)+\frac hR\eta'(h/R),
\]
and
\[
g_R''(h)
=
\frac2R\eta'(h/R)
+
\frac h{R^2}\eta''(h/R).
\]
There exists a constant $C>0$, independent
of $R$, such that
\[
|g_R'(h)|\leq C,
\qquad
|g_R''(h)|
\leq
\frac CR\mathbf1_{\{R<h<2R\}}.
\]
Moreover, both derivatives vanish for $h\leq R$.
Hence,
\[
L_{\phi,a}H_R
=
g_R'(h)L_{\phi,a}h
+
\phi g_R''(h)
\left(h-\frac1h\right).
\]
Since $h\leq2R$ on the support of
$g_R''$, it follows that
\[
|L_{\phi,a}H_R|
\leq
C(\Lambda+\mu)\mathbf1_{\{h>R\}},
\]
which proves~\eqref{energytailgenerator}.
The bounds~\eqref{energytailtestbounds}
follow directly from the definition of $\eta$.
\end{proof}

\begin{proposition}\label{energytaillemma}
Let $G$ be a weak solution on $[0,T]$ with initial
data $G_0$ satisfying $E_{G_0}<\infty$.
Then there exists a constant $C>0$, independent
of $G$, $T$, and $R$, such that
\begin{equation}\label{energytailestimate}
\sup_{0\leq t\leq T}
\mathcal T_{2R}(G(t))
\leq
\mathcal T_R(G_0)
+
\frac{C\mu ET}{R},
\qquad R>1,
\end{equation}
where $E=E_{G_0}$.
Consequently,
\begin{equation}\label{uniformenergytails}
\lim_{R\to\infty}
\sup_{0\leq t\leq T}
\mathcal T_R(G(t))=0.
\end{equation}
\end{proposition}

\begin{proof}
Let $H_R$ be the function introduced in
Lemma~\ref{energytailtestlemma}.
Since
\[
0<\Phi_G(t)\leq\frac{\mu}{3},
\qquad
|j_G(t)|<1,
\]
the generator estimate
\eqref{energytailgenerator}, applied with
$\Lambda=\mu/3$, gives
\[
|L_GH_R(p)|
\leq
C\mu\mathbf1_{\{\langle p\rangle>R\}}.
\]
We first justify the use of $H_R$ in the
integrated weak formulation.
Let $\chi\in C_c^\infty([0,\infty))$ satisfy
\[
0\leq\chi\leq1,
\qquad
\chi(r)=1\quad(r\leq1),
\qquad
\chi(r)=0\quad(r\geq2),
\]
and set
\[
\chi_N(p)=\chi\left(\frac{|p|}{N}\right).
\]
For fixed $R$ and $N$ sufficiently large,
\[
H_R(p)=\langle p\rangle
\]
on the support of $\nabla_p\chi_N$.
We compute
\[
\begin{aligned}
L_G(\chi_NH_R)
={}&
\chi_NL_GH_R
+
\Phi_GH_R\nabla_p\cdot(D\nabla_p\chi_N)\\
&+
2\Phi_G\nabla_p\chi_N\cdot D\nabla_pH_R
-
H_RW_G\cdot\nabla_p\chi_N.
\end{aligned}
\]
On the support of $\nabla_p\chi_N$, we have
\[
D\nabla_pH_R=p,
\qquad
|\nabla_p\chi_N|\leq\frac CN,
\qquad
\left|
\nabla_p\cdot(D\nabla_p\chi_N)
\right|\leq\frac CN.
\]
Moreover, \eqref{betanoen} implies
$|W_G(t,p)|\leq3\mu$.
Since $\langle p\rangle/N$ is bounded on
the support of $\nabla_p\chi_N$, it follows that
\[
\left|
L_G(\chi_NH_R)-\chi_NL_GH_R
\right|
\leq
C\mu\mathbf1_{\{N\leq|p|\leq2N\}}.
\]
Using $\chi_NH_R$ in
\eqref{integratedweak} and letting $N\to\infty$,
we obtain
\[
\int_{\R^3}H_R(p)G(t,p)\,dp
=
\int_{\R^3}H_R(p)G_0(p)\,dp+
\int_0^t\int_{\R^3}
L_GH_R(p)G(s,p)\,dp\,ds.
\]
Indeed, the cutoff errors vanish by dominated
convergence, using conservation of mass.
The other terms converge by dominated convergence,
since $H_R\leq\langle p\rangle$, the energy is
finite, and $L_GH_R$ is bounded.
By \eqref{energytailgenerator},
\[
\int_{\R^3}H_RG(t)\,dp\leq
\int_{\R^3}H_RG_0\,dp+
C\mu\int_0^t
\int_{\langle p\rangle>R}G(s,p)\,dp\,ds.
\]
Since the energy is conserved,
\[
\int_{\langle p\rangle>R}G(s,p)\,dp
\leq
\frac ER,
\]
and hence
\[
\int_{\R^3}H_RG(t)\,dp
\leq
\mathcal T_R(G_0)+\frac{C\mu Et}{R}.
\]
Finally, by \eqref{energytailtestbounds},
\[
\mathcal T_{2R}(G(t))
\leq
\int_{\R^3}H_RG(t)\,dp.
\]
Taking the supremum over $0\leq t\leq T$
proves \eqref{energytailestimate}.
Since $G_0$ has finite energy,
\[
\lim_{R\to\infty}\mathcal T_R(G_0)=0.
\]
Letting $R\to\infty$ in
\eqref{energytailestimate} therefore proves
\eqref{uniformenergytails}.
\end{proof}

%===========================================================
\subsection{Completion of the proof}
%===========================================================

We now complete the proof of Theorem~\ref{mainthm}.
We first show that the diffusion coefficient remains
uniformly positive on every finite time interval.
We then establish uniqueness and use the resulting
a priori estimates to continue the local solution
globally in time.

%===========================================================
\subsubsection*{Positivity of the diffusion coefficient}
%===========================================================

We start by showing that the nonlinear diffusion coefficient
cannot degenerate on a finite time interval.

Let $G$ be a weak solution of~\eqref{mainPDE} on $[0,T]$, where $T>0$.
By conservation of energy and momentum,
\[
\int_{\R^3}\langle p\rangle G(t,p)\,dp=E,
\qquad
\int_{\R^3}p\,G(t,p)\,dp=P,
\qquad 0\leq t\leq T.
\]
Moreover, by Lemma~\ref{Minkowskilemma},
\[
\Delta:=E^2-|P|^2-1>0.
\]
By Proposition~\ref{variancegap}, the inequality
$1+\gamma_E\leq2$, and the lower bound
\eqref{Phinoen}, for every $R>1$,
\[
\Phi_G(t)
\geq
\frac{9\mu}{128E^2R^2}
\left(
\Delta-4E\mathcal T_R(G(t))
\right)_+^2,
\qquad 0\leq t\leq T.
\]
which is equivalent to~ \eqref{mainPhibound}. By Proposition~\ref{energytaillemma}, we may choose
$R_T>1$, depending only on $G_0,E,P,\mu$ and $T$,
such that
\[
\sup_{0\leq t\leq T}
\mathcal T_{R_T}(G(t))
\leq\frac{\Delta}{8E}.
\]
Consequently,
\[
\Delta-4E\mathcal T_{R_T}(G(t))
\geq\frac{\Delta}{2},
\qquad 0\leq t\leq T.
\]
Hence
\begin{equation}\label{finitePhigap}
\Phi_G(t)
\geq
\frac{9\mu\Delta^2}{512E^2R_T^2}
=:\delta_T>0,
\qquad 0\leq t\leq T.
\end{equation}
Thus the nonlinear diffusion coefficient remains
uniformly positive on every finite time interval.

%===========================================================
\subsubsection*{Uniqueness}
%===========================================================

Next we prove uniqueness in the class of weak solutions
appearing in Theorem~\ref{mainthm}.

Let $G$ and $H$ be two weak solutions of~\eqref{mainPDE} on $[0,T]$
with the same initial data $G_0$.
By the conservation laws,
\[
\int_{\R^3}\langle p\rangle G(t,p)\,dp
=
\int_{\R^3}\langle p\rangle H(t,p)\,dp
=
E
\]
and
\[
\int_{\R^3}p\,G(t,p)\,dp
=
\int_{\R^3}p\,H(t,p)\,dp
=
P.
\]
Consequently, by Lemma~\ref{estimatesBeta},
\begin{equation}\label{uniquenessgamma}
|j_G(t)|,|j_H(t)|
\leq
\gamma_E
=
\sqrt{1-\frac1{E^2}}
<1.
\end{equation}
Since
\[
G,H\in C([0,T];L^1(\R^3)),
\]
the coefficient functions are continuous. Indeed,
\[
|j_G(t)-j_G(s)|
\leq
\|G(t)-G(s)\|_{L^1},
\]
and similarly for $H$. Moreover, by
\eqref{uniquenessgamma} and the definitions
of $\Phi_G$ and $\Phi_H$,
\[
|\Phi_G(t)-\Phi_G(s)|
\leq
C\|G(t)-G(s)\|_{L^1},
\]
and similarly for $\Phi_H$. 
By \eqref{finitePhigap}, there exists
$\delta_T>0$, depending only on $G_0,\mu$
and $T$, such that
\begin{equation}\label{uniquenessdelta}
\delta_T
\leq
\Phi_G(t),\Phi_H(t)
\leq
\frac{\mu}{3},
\qquad 0\leq t\leq T.
\end{equation}
Hence
\[
(\Phi_G,j_G),(\Phi_H,j_H)
\in
C([0,T];\R\times\R^3).
\]
Thus $G$ and $H$ may be regarded as solutions of the frozen
Fokker--Planck equations corresponding respectively to the
coefficient functions
$(\Phi_G,j_G)$ and $(\Phi_H,j_H)$.
These coefficient functions are continuous and satisfy the
bounds \eqref{uniquenessgamma}--\eqref{uniquenessdelta}.
Extending them constantly beyond $T$, we may apply
Theorem~\ref{frozenwellposedness} and
Lemma~\ref{frozenestimates} on $[0,T]$.

Let $P_t$ and $Q_t$ solve the stochastic differential equations
corresponding respectively to $(\Phi_G,j_G)$ and $(\Phi_H,j_H)$,
driven by the same Brownian motion, and with $P_0=Q_0$
almost surely. Let $G_0(p)\,dp$ be the probability distribution of $P_0$ and $Q_0$. 
By the uniqueness statement in
Theorem~\ref{frozenwellposedness}, which applies
to all weak solutions in the stated class,
the time marginal distributions of $P_t$ and $Q_t$
are respectively $G(t,p)\,dp$ and $H(t,p)\,dp$. 
Set
\[
Z_t=P_t-Q_t.
\]
The coefficient estimates in
Lemma~\ref{frozencoefflemma}, together with
$
\E\langle P_t\rangle
=
\E\langle Q_t\rangle
=
E,
$
give
\begin{equation}\label{uniquenessSDE}
\E|Z_t|^2
\leq
C_T\int_0^t\E|Z_s|^2\,ds+
C_T\int_0^t
\left(
|\Phi_G(s)-\Phi_H(s)|^2
+
|j_G(s)-j_H(s)|^2
\right)ds.
\end{equation}
Here the estimate is justified by localization,
as in the proof of Lemma~\ref{frozenestimates}.
In particular, it does not require finite second
moments of the individual processes.
Since the joint distribution of $(P_t,Q_t)$ is a coupling of
$G(t)$ and $H(t)$,
\begin{equation}\label{W1coupling}
W_1(G(t),H(t))
\leq
\E|P_t-Q_t|
\leq
\bigl(\E|Z_t|^2\bigr)^{1/2}.
\end{equation}
By Lemma~\ref{momentLipschitz} and
\eqref{uniquenessgamma},
\[
|\Phi_G(t)-\Phi_H(t)|
+
|j_G(t)-j_H(t)|
\leq
C_TW_1(G(t),H(t)).
\]
Hence
\[
|\Phi_G(t)-\Phi_H(t)|^2
+
|j_G(t)-j_H(t)|^2
\leq
C_TW_1(G(t),H(t))^2.
\]
Using \eqref{W1coupling}, we obtain from
\eqref{uniquenessSDE}
\[
\E|Z_t|^2
\leq
C_T\int_0^t\E|Z_s|^2\,ds.
\]
Since $Z_0=0$, Gronwall's inequality gives
$\E|Z_t|^2=0$ for $0\leq t\leq T$.
Therefore
\[
P_t=Q_t
\qquad\text{almost surely},
\]
and hence $G(t,p)=H(t,p)$
for almost every $p\in\R^3$ and every $t\in[0,T]$.

%===========================================================
\subsubsection*{Global continuation}
%===========================================================

We finally prove that the local solution in
Theorem~\ref{localexistencetheorem}
extends globally in time.

By local existence and uniqueness, the solution
can be continued to a unique maximal interval
$[0,T_*)$, where $0<T_*\leq\infty$.
Suppose, by contradiction, that $T_*<\infty$.
By conservation of energy,
\[
E_{G(t)}=E,
\qquad 0\leq t<T_*.
\]
Consequently, by Lemma~\ref{estimatesBeta},
\[
|j_G(t)|
\leq
\gamma_E
=
\sqrt{1-\frac1{E^2}}
<1,
\qquad 0\leq t<T_*.
\]
By Proposition~\ref{energytaillemma}, for every
$t<T_*$ and $R>1$,
\[
\mathcal T_{2R}(G(t))
\leq
\mathcal T_R(G_0)
+
\frac{C\mu ET_*}{R}.
\]
Since $G_0$ has finite energy, we may choose
$R_*>1$, depending only on
$G_0,E,P,\mu$ and $T_*$, such that
\[
\mathcal T_{R_*}(G(t))
\leq
\frac{\Delta}{8E},
\qquad 0\leq t<T_*,
\]
where
\[
\Delta=E^2-|P|^2-1>0.
\]
Proposition~\ref{variancegap} and
Lemma~\ref{Phiboundslemma} therefore give
\[
\Phi_G(t)
\geq
\delta_*
:=
\frac{9\mu\Delta^2}{512E^2R_*^2}
>0,
\qquad 0\leq t<T_*.
\]
Moreover,
\[
\Phi_G(t)\leq\frac{\mu}{3}.
\]
Fix $t_0<T_*$.
We restart the local fixed-point argument at time
$t_0$, taking $G(t_0)$ as initial data.
The estimates above give
\[
\Phi_{G(t_0)}\geq\delta_*>0,
\qquad
|j_{G(t_0)}|\leq\gamma_E<1,
\qquad
E_{G(t_0)}=E.
\]
By the uniform statement in
Theorem~\ref{localexistencetheorem},
there exists $\tau_*>0$, depending only on
$\delta_*,\gamma_E,E$ and $\mu$, and hence
independent of $t_0<T_*$, such that the local
existence argument with initial data $G(t_0)$
produces a solution on $[t_0,t_0+\tau_*]$.
Choose $t_0<T_*$ sufficiently close to $T_*$
that $t_0+\tau_*>T_*$.
By uniqueness, the restarted solution agrees
with $G$ on their common interval of existence.
By Lemma~\ref{integratedweaklemma}, the two solutions
can be concatenated to give a weak solution extending
$G$ beyond $T_*$, contradicting the maximality of
$T_*$.
Hence $T_*=\infty$.

This proves global existence and completes the proof of
Theorem~\ref{mainthm}.

%===========================================================
\section{Proof of Theorem~\ref{regularizationthm}}
%===========================================================

We now prove Theorem~\ref{regularizationthm}.
Fix $T>0$ and set
\[
\phi(t)=\Phi_G(t),
\qquad
a(t)=j_G(t).
\]
By Theorem~\ref{mainthm}, there exists a constant
$\delta_T>0$ such that
\[
0<\delta_T\leq\phi(t)\leq\frac{\mu}{3},
\qquad
|a(t)|\leq\gamma_E,
\qquad
0\leq t\leq T,
\]
where
\[
\gamma_E=\sqrt{1-\frac1{E^2}}<1.
\]
The diffusion matrix of the Fokker--Planck
equation is $\phi(t)D(p)$.
By~\eqref{localellipticity}, for
$|p|\leq R$ and $\xi\in\R^3$,
\begin{equation}\label{localellipticity2}
\frac{\delta_T}{\sqrt{1+R^2}}|\xi|^2
\leq
\phi(t)\xi^TD(p)\xi
\leq
\frac{\mu\sqrt{1+R^2}}{3}|\xi|^2.
\end{equation}
Thus, the PDE~\eqref{mainPDE} is uniformly parabolic
on every bounded momentum set, uniformly on finite
time intervals.

We divide the proof of
Theorem~\ref{regularizationthm} into several subsections.

%===========================================================
\subsection{An algebra of admissible observables}
%===========================================================

Fix $0<\gamma<1$. For $|a|\leq\gamma$, set
\[
z=1-a\cdot v.
\]
Then
\[
z\geq1-\gamma>0.
\]
For $k\in\mathbb N_0$, let $\mathcal A_k$ be the finite linear span
of functions of the form
\[
f(a,p)
=
\langle p\rangle^{-k}
\frac{P(a,v)}{z^\ell},
\qquad
\ell\in\mathbb N_0,
\]
where $P$ is a polynomial in $(a,v)$. We call finite sums of elements
of the classes $\mathcal A_k$ admissible observables.
Every function in $\mathcal A_k$, together with all its derivatives
with respect to $a$, is bounded uniformly on
$\{|a|\leq\gamma\}\times\R^3$.
Moreover,
\begin{equation}\label{Aproduct}
\mathcal A_k\mathcal A_m
\subset\mathcal A_{k+m},
\qquad k,m\in\mathbb N_0.
\end{equation}
Since $\langle p\rangle^{-2}=1-|v|^2$ 
is a polynomial in $v$, we also have
\begin{equation}\label{Anesting}
\mathcal A_{k+2}\subset\mathcal A_k,
\qquad k\in\mathbb N_0.
\end{equation}
The identities
\[
\partial_{p_j}\langle p\rangle^{-k}
=
-k\langle p\rangle^{-k-1}v_j,
\quad
\partial_{p_j}v_i
=
\langle p\rangle^{-1}
(\delta_{ij}-v_iv_j),
\]
and
\[
\partial_{p_j}z
=
-\langle p\rangle^{-1}
\bigl(a_j-(a\cdot v)v_j\bigr)
\]
give
\begin{equation}\label{Aderivatives}
\partial_{p_j}\mathcal A_k
\subset\mathcal A_{k+1},
\qquad
\partial_{a_j}\mathcal A_k
\subset\mathcal A_k.
\end{equation}

We shall also use observables depending polynomially on an additional
scalar parameter $\phi$. More precisely, let $\mathcal A$ denote the
finite linear span of functions of the form
\[
f(\phi,a,p)=\phi^r g(a,p),
\qquad
r\in\mathbb N_0,
\qquad
g\in\mathcal A_k
\]
for some $k\geq0$. We shall refer to the elements of $\mathcal A$
as extended admissible observables. By \eqref{Aproduct},
$\mathcal A$ is an algebra under pointwise multiplication.
On sets where $|a|\leq\gamma<1$ and $\phi$ remains bounded,
every extended admissible observable, together with all its
derivatives with respect to $\phi$ and $a$, is uniformly bounded.

Recall that the generator of the frozen Fokker--Planck
equation is
\[
L_{\phi,a}f
=
\phi\nabla_p\cdot(D\nabla_pf)
-
V_a\cdot\nabla_pf,
\]
where
\[
D(p)=\frac{I+p\otimes p}{\langle p\rangle},
\qquad
V_a(p)
=
\mu\left(
v-\frac{a}{\langle p\rangle^2z}
\right),
\qquad
z=1-a\cdot v.
\]

\begin{lemma}\label{closurelemma}
For every $k\geq0$ and each fixed $\phi$,
\begin{equation}\label{generatorclosure}
L_{\phi,a}\mathcal A_k
\subset\mathcal A_{k+1}.
\end{equation}
The coefficients appearing in this inclusion depend
polynomially on $\phi$. In particular, the class
$\mathcal A$ is preserved by $L_{\phi,a}$,
$\partial_\phi$, and $\nabla_a$.
\end{lemma}

\begin{proof}
Let $f\in\mathcal A_k$. By \eqref{Aderivatives},
\[
\partial_{p_i}f\in\mathcal A_{k+1},
\qquad
\partial_{p_i}\partial_{p_j}f
\in\mathcal A_{k+2}.
\]
Since $D_{ij}=\langle p\rangle^{-1}\delta_{ij}+\langle p\rangle v_iv_j$,
we have
$D_{ij}\partial_{p_i}\partial_{p_j}f\in\mathcal A_{k+3}+\mathcal A_{k+1}$, and thus, by \eqref{Anesting},
\[
D_{ij}\partial_{p_i}\partial_{p_j}f
\in\mathcal A_{k+1}.
\]
Moreover, since
$\partial_{p_i}D_{ij}=3v_j$, 
we have
$(\partial_{p_i}D_{ij})\partial_{p_j}f
\in\mathcal A_{k+1}$.
It follows that
\[
\nabla_p\cdot(D\nabla_pf)
\in\mathcal A_{k+1}.
\]
Next, using
\[
V_a
=
\mu\left(
v-\frac{a}{\langle p\rangle^2z}
\right),
\]
we obtain
\[
V_a\cdot\nabla_pf
\in
\mathcal A_{k+1}+\mathcal A_{k+3}
\subset\mathcal A_{k+1}.
\]
Combining these inclusions gives
\[
L_{\phi,a}f\in\mathcal A_{k+1},
\]
which proves \eqref{generatorclosure}.
Finally, the operator $L_{\phi,a}$ depends linearly
on $\phi$, while its remaining coefficients belong
to the classes of admissible observables. Therefore
$L_{\phi,a}$ maps extended admissible observables
into extended admissible observables.
The invariance of $\mathcal A$ under
$\partial_\phi$ and $\nabla_a$ follows directly
from its definition and \eqref{Aderivatives}.
\end{proof}

%===========================================================
\subsection{First derivatives of the nonlinear coefficients}
%===========================================================

We first establish differentiability of $a$ and $\phi$ directly from
the weak equation.

Recall that
\[
a(t)
=j_G(t)=
\int_{\R^3}v(p)G(t,p)\,dp.
\]
A direct calculation gives
\[
D\nabla_pv_i
=
\frac{e_i}{\langle p\rangle^2},
\qquad
\nabla_p\cdot(D\nabla_pv_i)
=
-\frac{2v_i}{\langle p\rangle^3},
\]
and
\[
V_a\cdot\nabla_pv_i
=
\frac{\mu(v_i-a_i)}
{\langle p\rangle^3(1-a\cdot v)}.
\]
Using suitable cutoff approximations of $v_i$
in the integrated weak formulation, we obtain
\begin{equation}\label{aprimefinal}
a'(t)
=
-2\phi(t)
\int_{\R^3}
\frac{v}{\langle p\rangle^3}G(t,p)\,dp
-
\mu
\int_{\R^3}
\frac{v-a(t)}
{\langle p\rangle^3(1-a(t)\cdot v)}
G(t,p)\,dp.
\end{equation}
By conservation of energy and
Lemma~\ref{estimatesBeta},
$|a(t)|\leq\gamma_E<1,
\qquad t\geq0$. 
The integrands in \eqref{aprimefinal} are therefore
bounded continuous functions of $(a,p)$, uniformly
for $|a|\leq\gamma_E$.
Since
$G\in C([0,\infty);L^1(\R^3))$
and $a,\phi$ are continuous, the right hand side of
\eqref{aprimefinal} is continuous. Hence,
\[
a\in C^1((0,\infty)).
\]
Next set
\[
F_a(p)
=
\frac{1}
{\langle p\rangle^2(1-a\cdot v)}.
\]
Then $F_a\in\mathcal A_2$ and
\[
\phi(t)
=
\frac{\mu}{3}
\left[
1-
\int_{\R^3}F_{a(t)}(p)G(t,p)\,dp
\right].
\]
Since $a\in C^1((0,\infty))$, we may use cutoff
approximations to the time-dependent observable
$F_{a(t)}$ in the integrated weak formulation.
The chain rule in the parameter $a$ gives
\begin{equation}\label{phiprimefinal}
\phi'(t)
=
-\frac{\mu}{3}
\left[
\int_{\R^3}
L_{\phi(t),a(t)}F_{a(t)}(p)
G(t,p)\,dp
+
a'(t)\cdot
\int_{\R^3}
\nabla_aF_{a(t)}(p)
G(t,p)\,dp
\right],
\end{equation}
where
\[
\nabla_aF_a
=
\frac{v}
{\langle p\rangle^2(1-a\cdot v)^2}.
\]
By Lemma~\ref{closurelemma},
\[
L_{\phi,a}F_a\in\mathcal A_3.
\]
The functions $L_{\phi,a}F_a$ and $\nabla_aF_a$
are therefore uniformly bounded and continuous
in their parameters on
\[
\left\{
0\leq\phi\leq\frac{\mu}{3},
\quad |a|\leq\gamma_E
\right\}.
\]
It follows that the right hand side of
\eqref{phiprimefinal} is continuous. Hence,
\[
\phi\in C^1((0,\infty)).
\]

%===========================================================
\subsection{Differentiation of admissible moments}
%===========================================================

For an extended admissible observable $f=f(\phi,a,p)$, define
\[
\mathcal M_f(t)
=
\int_{\R^3}
f(\phi(t),a(t),p)G(t,p)\,dp.
\]
\begin{lemma}\label{momentdifferentiation}
Suppose that $a,\phi\in C^1(I)$ on a compact interval
$I\subset(0,\infty)$ and that $f$ is an extended admissible
observable. Then $\mathcal M_f\in C^1(I)$ and
\begin{align}
\frac d{dt}\mathcal M_f(t)
={}&
\int_{\R^3}
L_{\phi(t),a(t)}
f(\phi(t),a(t),p)
G(t,p)\,dp
\nonumber\\
&+
\phi'(t)
\int_{\R^3}
\partial_\phi f(\phi(t),a(t),p)
G(t,p)\,dp
\nonumber\\
&+
a'(t)\cdot
\int_{\R^3}
\nabla_af(\phi(t),a(t),p)
G(t,p)\,dp.
\label{expectationrule}
\end{align}
\end{lemma}

\begin{proof}
Let
\[
\chi\in C_c^\infty([0,\infty)),
\qquad
0\leq\chi\leq1,
\]
satisfy
\[
\chi(r)=1\quad\text{for }0\leq r\leq1,
\qquad
\chi(r)=0\quad\text{for }r\geq2,
\]
and set
\[
\chi_R(p)=\chi\left(\frac{|p|}{R}\right),
\qquad R>1.
\]
Then
\[
|\nabla_p\chi_R|\leq\frac{C}{R},
\]
and all derivatives of $\chi_R$ are supported in
\[
A_R=\{R\leq|p|\leq2R\}.
\]
Since $\chi_R$ is radial and the radial eigenvalue of $D(p)$ is
$\langle p\rangle$,
\[
D(p)\nabla_p\chi_R
=
\langle p\rangle\nabla_p\chi_R.
\]
A direct calculation therefore gives
\begin{equation}\label{cutoffbounds}
|D\nabla_p\chi_R|\leq C,
\qquad
\left|
\nabla_p\cdot(D\nabla_p\chi_R)
\right|
\leq\frac{C}{R}.
\end{equation}

Suppose first that $f\in\mathcal A_k$. Uniformly for
$|a|\leq\gamma<1$,
\[
|f(a,p)|\leq C\langle p\rangle^{-k},
\qquad
|\nabla_pf(a,p)|
\leq C\langle p\rangle^{-k-1}.
\]
Since the largest eigenvalue of $D(p)$ is $\langle p\rangle$,
\[
|D\nabla_pf|
\leq C\langle p\rangle^{-k}.
\]
Hence, 
\begin{align}
L_{\phi,a}(\chi_Rf)-\chi_RL_{\phi,a}f
={}&
\phi f\,\nabla_p\cdot(D\nabla_p\chi_R)
\nonumber\\
&+
2\phi\nabla_p\chi_R\cdot D\nabla_pf
-
fV_a\cdot\nabla_p\chi_R.
\label{cutoffcommutator}
\end{align}
Using the estimates
\[
|a(t)|\leq\gamma_E<1,
\qquad
0<\phi(t)\leq\frac{\mu}{3}
\]
in Lemmas~\ref{estimatesBeta}
and~\ref{Phiboundslemma}, we obtain $
|V_{a(t)}(p)|\leq C_E$. 
It follows from \eqref{cutoffbounds} that
\[
\left|
L_{\phi,a}(\chi_Rf)-\chi_RL_{\phi,a}f
\right|
\leq
\frac{C_I}{R}
\langle p\rangle^{-k}\mathbf1_{A_R}(p).
\]
Hence
\[
\int_{\R^3}
\left|
L_{\phi,a}(\chi_Rf)-\chi_RL_{\phi,a}f
\right|
G(t,p)\,dp
\leq
\frac{C_I}{R},
\]
uniformly for $t\in I$.
By linearity, the same estimate holds for
extended admissible observables, with a constant
depending also on the observable $f$.
Thus, all cutoff commutators vanish as $R\to\infty$.
Since $a,\phi\in C^1(I)$, the function
\[
(t,p)\longmapsto
\chi_R(p)f(\phi(t),a(t),p)
\]
belongs to $C^{1,2}(I\times\R^3)$ and has compact
momentum support uniformly on $I$.
It is therefore an admissible test function in the
integrated weak formulation of
Lemma~\ref{integratedweaklemma}.
We compute
\[
\partial_t f
=
\phi'(t)\partial_\phi f
+
a'(t)\cdot\nabla_af.
\]

For $t_1,t_2\in I$, with $t_1<t_2$,
the integrated weak formulation yields
\[
\begin{aligned}
&\int_{\R^3}
\chi_R f(\phi(t_2),a(t_2),p)
G(t_2,p)\,dp\\
&\quad-
\int_{\R^3}
\chi_R f(\phi(t_1),a(t_1),p)
G(t_1,p)\,dp\\
&=
\int_{t_1}^{t_2}\int_{\R^3}
\chi_R
\left[
L_{\phi,a}f
+
\phi'\partial_\phi f
+
a'\cdot\nabla_af
\right]G\,dp\,dt\\
&\quad+
\int_{t_1}^{t_2}\int_{\R^3}
\left[
L_{\phi,a}(\chi_Rf)
-
\chi_RL_{\phi,a}f
\right]G\,dp\,dt,
\end{aligned}
\]
where the coefficients and observables in the
time integrals are evaluated at $(t,p)$.
The last term tends to zero as $R\to\infty$
by the commutator estimate.
Moreover, Lemma~\ref{closurelemma} shows that
\[
L_{\phi,a}f,\qquad
\partial_\phi f,\qquad
\nabla_af
\]
are extended admissible observables and hence
are uniformly bounded on the relevant parameter set.
Since $a',\phi'$ are bounded on $I$,
dominated convergence allows us to pass to the
limit in the remaining terms.
We obtain
\[
\mathcal M_f(t_2)-\mathcal M_f(t_1)
=
\int_{t_1}^{t_2}B_f(t)\,dt,
\]
where
\[
\begin{aligned}
B_f(t)
={}&
\int_{\R^3}
L_{\phi(t),a(t)}
f(\phi(t),a(t),p)
G(t,p)\,dp\\
&+
\phi'(t)
\int_{\R^3}
\partial_\phi f(\phi(t),a(t),p)
G(t,p)\,dp\\
&+
a'(t)\cdot
\int_{\R^3}
\nabla_af(\phi(t),a(t),p)
G(t,p)\,dp.
\end{aligned}
\]
The function $B_f$ is continuous on $I$,
because all the observables involved are
bounded and continuous in their parameters, while
\[
G\in C(I;L^1(\R^3)),
\qquad
a,\phi\in C^1(I).
\]
Therefore $\mathcal M_f\in C^1(I)$ and
\[
\mathcal M_f'(t)=B_f(t),
\]
which proves \eqref{expectationrule}.
\end{proof}

%===========================================================
\subsection{Smoothness of the nonlinear coefficients}
%===========================================================

\begin{proposition}\label{coeffsmoothprop}
The nonlinear coefficient functions satisfy
\begin{equation}\label{coeffsmooth}
a,\phi\in C^\infty((0,\infty)).
\end{equation}
\end{proposition}

\begin{proof}
We have already proved that
\[
a,\phi\in C^1((0,\infty)).
\]
Fix a compact interval $I\subset(0,\infty)$.
By conservation of energy and
Lemma~\ref{estimatesBeta}, $
|a(t)|\leq\gamma_E<1$. 
We prove by induction on $m\geq1$ that
\[
a,\phi\in C^m(I),
\qquad
\mathcal M_f\in C^m(I)
\]
for every extended admissible observable $f$.
The case $m=1$ follows from the preceding
subsections and Lemma~\ref{momentdifferentiation}.
Suppose that the assertion holds for some
$m\geq1$. By \eqref{aprimefinal}, the function
$a'$ is a finite linear combination of extended
admissible moments. Hence
\[
a\in C^{m+1}(I).
\]
Using \eqref{phiprimefinal}, the regularity of
$a'$, and the induction hypothesis for
admissible moments, we similarly obtain
\[
\phi'\in C^m(I),
\qquad
\phi\in C^{m+1}(I).
\]

It remains to improve the regularity of the
admissible moments. By Lemma~\ref{momentdifferentiation},
\[
\frac d{dt}\mathcal M_f
=
\mathcal M_{L_{\phi,a}f}
+
\phi'\mathcal M_{\partial_\phi f}
+
a'\cdot\mathcal M_{\nabla_a f}.
\]
By Lemma~\ref{closurelemma}, all the
observables appearing on the right hand side
are extended admissible. Their moments belong
to $C^m(I)$ by the induction hypothesis.
Since $a',\phi'\in C^m(I)$, it follows that
\[
\mathcal M_f\in C^{m+1}(I).
\]
This completes the induction.
Therefore
\[
a,\phi\in C^\infty(I).
\]
Since $I\subset(0,\infty)$ was arbitrary,
the result follows.
\end{proof}

%===========================================================
\subsection{Interior smoothness}
%===========================================================

We now prove that the weak solution is smooth for positive times.
Fix
\[
0<\tau<T,
\qquad
0<R<R',
\]
and set
\[
Q'=(\tau/2,T)\times B_{R'}.
\]
On $Q'$, the function $G$ is a distributional solution of
\begin{equation}\label{localPDEreg}
\partial_tG
=
\partial_i
\left(
\phi(t)D^{ij}(p)\partial_jG
+
V_{a(t)}^i(p)G
\right).
\end{equation}
Set
\[
A(t,p)=\phi(t)D(p).
\]
By \eqref{localellipticity2}, the matrix $A$
is uniformly positive definite on $Q'$.
Moreover, Proposition~\ref{coeffsmoothprop}
implies that all the coefficients in
\eqref{localPDEreg} are $C^\infty$ on $Q'$.
Let
\[
S(t,p)=A(t,p)^{1/2}.
\]
Since $A$ is smooth and positive definite on $Q'$,
its positive square root $S$ is smooth.
Define the vector fields
\[
X_k
=
\sum_{i=1}^3S_{ik}(t,p)\partial_{p_i},
\qquad k=1,2,3.
\]
Since $SS^T=A$, we have
\[
\sum_{k=1}^3X_k^2
=
A^{ij}\partial_{p_i}\partial_{p_j}
+
b^j(t,p)\partial_{p_j},
\]
where
\[
b^j(t,p)
=
\sum_{i,k=1}^3
S_{ik}(t,p)\partial_{p_i}S_{jk}(t,p)
\]
are smooth functions on $Q'$.
Expanding \eqref{localPDEreg}, we see that
$G$ satisfies
\[
PG=0
\]
in the sense of distributions on $Q'$, where
\[
P
=
-\partial_t
+
A^{ij}\partial_{p_i}\partial_{p_j}
+
\left(
\partial_{p_i}A^{ij}+V_{a(t)}^j
\right)\partial_{p_j}
+
\partial_{p_i}V_{a(t)}^i.
\]
Consequently,
\[
P
=
\sum_{k=1}^3X_k^2+X_0+c,
\]
where
\[
X_0=-\partial_t+Y,
\]
with
\[
Y
=
\sum_{j=1}^3
\left(
\partial_{p_i}A^{ij}
+
V_{a(t)}^j
-
b^j
\right)\partial_{p_j},
\]
and
\[
c=\partial_{p_i}V_{a(t)}^i.
\]
Thus $X_0,X_1,X_2,X_3$ are smooth vector fields
on $Q'$, and $c$ is a smooth function.
Since $A$ is positive definite, $S$ is invertible.
Therefore
\[
\operatorname{span}\{X_1,X_2,X_3\}
=
\operatorname{span}
\{\partial_{p_1},\partial_{p_2},\partial_{p_3}\}.
\]
Since $X_0$ has a nonzero component in the time direction,
it follows that
\[
\operatorname{span}\{X_0,X_1,X_2,X_3\}
=
T_{(t,p)}(\R\times\R^3)
\]
at every $(t,p)\in Q'$.
Thus, H\"ormander's condition is satisfied.
Finally,
\[
G\in C([0,\infty);L^1(\R^3))
\]
implies $G\in L^1_{\rm loc}(Q')$, and hence $G$ defines a
distribution on $Q'$. H\"ormander's hypoellipticity theorem
\cite[Theorem~1.1]{Hormander1967} therefore implies that
$G$ is smooth in the interior of $Q'$.
In particular,
\[
G\in C^\infty((\tau,T)\times B_R).
\]
Since $\tau,T,R$ are arbitrary, we conclude that
\begin{equation}\label{smoothfinal}
G\in C^\infty((0,\infty)\times\R^3).
\end{equation}

%===========================================================
\subsection{Strict positivity}
%===========================================================

We prove that $G$ is strictly positive by using
the parabolic strong maximum principle.

Fix $t>0$ and $p_*\in\R^3$, and set $s=t/2$.
Since $G$ has unit mass and is smooth by
\eqref{smoothfinal}, there exists $p_0\in\R^3$
such that
\[
G(s,p_0)>0.
\]
Choose $R>\max\{|p_0|,|p_*|\}$ and set
$\Omega=B_R$.
On $[s/2,2t]\times\overline{\Omega}$,
the equation is uniformly parabolic with
smooth coefficients.
Writing
\[
\phi(\tau)=\Phi_G(\tau),
\qquad
a(\tau)=j_G(\tau),
\]
the equation takes the form
\[
\partial_\tau G
=
\phi D:\nabla_p^2G
+
(3\phi v+V_a)\cdot\nabla_pG
+
(\nabla_p\cdot V_a)G.
\]
Since the coefficients are smooth, we may
choose $C>0$ such that
\[
\nabla_p\cdot V_{a(\tau)}(p)-C<0
\]
on $[s/2,2t]\times\overline{\Omega}$.
Define
\[
H(\tau,p)=e^{-C\tau}G(\tau,p).
\]
Then $H\geq0$ and
\[
\partial_\tau H
=
\phi D:\nabla_p^2H
+
(3\phi v+V_a)\cdot\nabla_pH
+
(\nabla_p\cdot V_a-C)H.
\]
Suppose that $G(t,p_*)=0$.
Then $H(t,p_*)=0$, so $-H$ attains
its maximum value zero at $(t,p_*)$.
By the parabolic strong maximum principle
\cite[Chapter~II, Theorem~2.7]{Lieberman},
$H$ must vanish throughout
$(s/2,t]\times\Omega$.
In particular, $H(s,p_0)=0$, contradicting
$G(s,p_0)>0$.
Therefore $G(t,p_*)>0$.
Since $t>0$ and $p_*\in\R^3$ were arbitrary,
we conclude that
\[
G(t,p)>0,
\qquad t>0,\quad p\in\R^3.
\]
This completes the proof of
Theorem~\ref{regularizationthm}.

%===========================================================
\section{Proof of Theorem~\ref{uniformPhitheorem}}
\label{longtimesec}
%===========================================================

In this section we prove Theorem~\ref{uniformPhitheorem}.
The main ingredient is a comparison principle for the
particle energy distributions of even solutions.

Assume that the initial data is an even probability density. Since the PDE~\eqref{mainPDE} is invariant under $p\mapsto-p$, uniqueness implies that
\[
G(t,p)=G(t,-p),\quad\text{for all $t>0$}
\]
and therefore
\[
j_G(t)=0,\quad\text{for all $t\geq0$.}
\]
For an even solution $G$, define
\[
\rho_G(t,z)
=
z\sqrt{z^2-1}
\int_{\mathbb S^2}
G\bigl(t,\sqrt{z^2-1}\,\omega\bigr)\,d\omega,
\qquad z>1.
\]
The function $\rho_G$ is the probability density
of the particle energy $z=\langle p\rangle$.
Indeed, passing to spherical coordinates and using
$|p|=\sqrt{z^2-1}$, we obtain
\[
\int_{\R^3}f(\langle p\rangle)G(t,p)\,dp
=
\int_1^\infty f(z)\rho_G(t,z)\,dz
\]
for every smooth function $f$ such that the integral in the left hand side is finite.
In particular, taking $f=1$ and $f=z$ gives
\begin{equation}\label{rhoG}
\int_1^\infty\rho_G(t,z)\,dz=1,\quad\int_1^\infty z\rho_G(t,z)\,dz=E.
\end{equation}
The tail distribution of $\rho_G$ is the decreasing function
\[
F_G(t,z)=\int_z^\infty\rho_G(t,s)\,ds,
\qquad z\geq1.
\]
By~\eqref{rhoG} and Tonelli's theorem,
\begin{equation}\label{energycumulativemoments}
F_G(t,1)=1,
\qquad
\int_1^\infty F_G(t,z)\,dz=E-1.
\end{equation}
Consequently,
\[
\Phi_G(t)
=
\frac{\mu}{3}
\int_1^\infty(1-z^{-2})\rho_G(t,z)\,dz.
\]
Integration by parts gives
\begin{equation}\label{PhiFcumulative}
\Phi_G(t)
=
\frac{2\mu}{3}
\int_1^\infty z^{-3}F_G(t,z)\,dz.
\end{equation}
A direct calculation using~\eqref{mainPDE} shows
that $\rho_G$ satisfies
\begin{equation}\label{energyFPnew}
\partial_t\rho_G
=
\Phi_G(t)\partial_z^2
\left[(z-z^{-1})\rho_G\right]
-
\partial_z
\left[
\left(3\Phi_G(t)-\mu(1-z^{-2})\right)
\rho_G
\right].
\end{equation}
Integrating this equation from $z$ to infinity
yields
\begin{equation}\label{energyFnew}
\partial_tF_G
=
\mathcal L_{\Phi_G(t)}F_G,
\end{equation}
where
\begin{equation}\label{energyLnew}
\mathcal L_\phi f
=
\phi(z-z^{-1})f''
+
\left[
\mu(1-z^{-2})-\phi(2-z^{-2})
\right]f'.
\end{equation}
These equations hold without assuming spherical
symmetry of $G$.

We first establish the uniform integrability
of the energy tails when $1<E<2$.

\begin{lemma}\label{lowenergyPDElemma}
Let $G_0$ be an even probability density with
finite energy satisfying $1<E<2$, and let $G$
be the global solution of~\eqref{mainPDE} with initial data $G_0$. Then
\begin{equation}\label{lowenergyuniformtails}
\lim_{R\to\infty}
\sup_{t\geq0}\mathcal T_R(G(t))=0.
\end{equation}
Consequently,
\[
\inf_{t\geq0}\Phi_G(t)>0.
\]
\end{lemma}

\begin{proof}
By~\eqref{elementaryPhi} and conservation
of energy,
\[
0<\Phi_G(t)\leq\phi_*,
\qquad
\phi_*=\frac{\mu}{3}(1-E^{-2}).
\]
Set
\[
\alpha=\frac{\mu}{\phi_*}
=\frac{3}{1-E^{-2}}.
\]
Since $1<E<2$, we have $\alpha>4$.
Define
\[
\pi(z)=z^{1-\alpha}\sqrt{z^2-1},
\qquad
S(z)=\int_z^\infty\pi(s)\,ds.
\]
A direct calculation using~\eqref{energyLnew}
gives
$\mathcal L_{\phi_*}S=0$. 
Moreover, for every $0<\phi\leq\phi_*$,
\[
\mathcal L_\phi S
=
\frac{\phi}{\phi_*}\mathcal L_{\phi_*}S
+
\mu\left(1-\frac{\phi}{\phi_*}\right)
(1-z^{-2})S'(z)
=
-\mu\left(1-\frac{\phi}{\phi_*}\right)
(1-z^{-2})\pi(z)
\leq0.
\]
Fix $\tau>0$ and write
\[
\rho_\tau(z)=\rho_G(\tau,z),
\qquad
F_\tau(z)=F_G(\tau,z).
\]
Define
\[
A_\tau(z)
=
\int_1^z\frac{\rho_\tau(s)}{S(s)}\,ds,
\qquad
\overline F_\tau(z)
=
F_\tau(z)+S(z)A_\tau(z).
\]
Since $\tau>0$, and by the regularity of $G$ proved in Theorem~\ref{regularizationthm}, the functions
$\rho_\tau$ and $F_\tau$ are smooth
on $(1,\infty)$, and consequently
$\overline F_\tau\in C^\infty((1,\infty))$.
By construction,
\[
\overline F_\tau(1)=1,
\qquad
\overline F_\tau(z)\geq F_\tau(z)\geq0.
\]
Differentiating, we obtain
\[
\overline F_\tau'(z)
=
-\pi(z)A_\tau(z)\leq0
\]
and
\[
\mathcal L_\phi\overline F_\tau
=
A_\tau\mathcal L_\phi S
-
\phi(z-z^{-1})
\frac{\pi(z)}{S(z)}\rho_\tau(z)
\leq0
\]
for every $0<\phi\leq\phi_*$.
Thus, $\overline F_\tau$ is a smooth
stationary supersolution of
\eqref{energyFnew} on $(1,\infty)$.
We now apply the classical parabolic
comparison principle. Fix $T>\tau$
and $1+\varepsilon<L$, and set
\[
w(t,z)=F_G(t,z)-\overline F_\tau(z).
\]
Since $\Phi_G(t)\leq\phi_*$, we have
\[
\partial_tw-\mathcal L_{\Phi_G(t)}w
=
\mathcal L_{\Phi_G(t)}\overline F_\tau
\leq0.
\]
By Theorem~\ref{mainthm}, $\Phi_G$ is
bounded away from zero on $[0,T]$, hence the equation is uniformly parabolic
on $[\tau,T]\times[1+\varepsilon,L]$.
At the initial time, $w(\tau,z)\leq0$.
At the left boundary, since $F_G(t,1)=1$,
\[
w(t,1+\varepsilon)
\leq1-\overline F_\tau(1+\varepsilon).
\]
At the right boundary, conservation of
energy gives
\[
w(t,L)\leq F_G(t,L)\leq\frac EL.
\]
The classical parabolic comparison principle
implies
\[
w(t,z)
\leq
\max\left\{
1-\overline F_\tau(1+\varepsilon),
\frac EL
\right\}
\]
for $\tau\leq t\leq T$ and
$1+\varepsilon\leq z\leq L$.
Since $\overline F_\tau$ is continuous
at $z=1$ and $\overline F_\tau(1)=1$,
letting $\varepsilon\to0$ and
$L\to\infty$ yields
\[
F_G(t,z)\leq\overline F_\tau(z),
\qquad t\geq\tau,\quad z\geq1.
\]
We now let $\tau\to0$.
Since $
G(\tau)\rightarrow G_0$ in $L^1(\R^3)$, 
we have
\[
\rho_\tau\rightarrow\rho_0
\quad\text{in }L^1((1,\infty)),
\]
where $\rho_0$ is the particle energy
probability density associated with $G_0$.
In particular,
\[
F_\tau(z)\rightarrow F_0(z)
\]
for every $z\geq1$.
Moreover, since $S$ is strictly
positive on $[1,\infty)$, for every
fixed $z>1$,
\[
A_\tau(z)
\rightarrow
A_0(z)
:=
\int_1^z\frac{\rho_0(s)}{S(s)}\,ds.
\]
Consequently,
\[
\overline F_\tau(z)
\rightarrow
\overline F(z)
:=
F_0(z)+S(z)A_0(z).
\]
We conclude that
\begin{equation}\label{lowenergyFcomparison}
F_G(t,z)\leq\overline F(z),
\qquad t\geq0,\quad z\geq1.
\end{equation}
It remains to prove that $\overline F$
has uniformly integrable energy tails.
Define
\[
\overline\rho(z)
=
-\overline F'(z)
=
\pi(z)\int_1^z
\frac{\rho_0(s)}{S(s)}\,ds.
\]
By Tonelli's theorem,
\[
\int_1^\infty\overline\rho(z)\,dz
=
\int_1^\infty
\frac{\rho_0(s)}{S(s)}
\left(\int_s^\infty\pi(z)\,dz\right)ds=1.
\]
Since $\overline F(1)=1$, this implies
\[
\overline F(z)
=
\int_z^\infty\overline\rho(s)\,ds.
\]
Furthermore,
\[
\int_1^\infty z\overline\rho(z)\,dz
=
\int_1^\infty
\frac{\rho_0(s)}{S(s)}
\Big(\int_s^\infty z\pi(z)\,dz\Big)\,ds.
\]
Since $\pi(z)\sim z^{2-\alpha}$ as $z\to\infty$
and $\alpha>4$, there exists a constant
$C_\alpha<\infty$ such that
\[
\frac{1}{S(s)}\int_s^\infty z\pi(z)\,dz\leq C_\alpha s,
\qquad s\geq1.
\]
Consequently,
\[
\int_1^\infty z\overline\rho(z)\,dz
\leq
C_\alpha\int_1^\infty
s\rho_0(s)\,ds
=
C_\alpha E<\infty.
\]
Finally, by~\eqref{lowenergyFcomparison},
\[
\mathcal T_R(G(t))=
RF_G(t,R)
+
\int_R^\infty F_G(t,z)\,dz\leq
R\overline F(R)
+
\int_R^\infty\overline F(z)\,dz=
\int_R^\infty z\overline\rho(z)\,dz.
\]
The right hand side is independent of
time and tends to zero as $R\to\infty$.
Therefore
\[
\lim_{R\to\infty}
\sup_{t\geq0}\mathcal T_R(G(t))=0.
\]
The uniform positivity of $\Phi_G$
now follows from
Corollary~\ref{uniformtailPhigap}.
\end{proof}

We next show that radial monotonicity is
preserved by the equation.

\begin{lemma}\label{radialmonotonicitylemma}
Let $H_0\in L^1(\R^3)$ be a radially
nonincreasing probability density with finite
energy, and let $H$ be the corresponding solution of~\eqref{mainPDE}. Then
\[
H(t,p)=h(t,|p|),
\]
where
\[
\partial_rh(t,r)\leq0,
\qquad t>0,\ r>0.
\]
\end{lemma}

\begin{proof}
Rotational invariance and uniqueness imply that
$H(t,p)=h(t,r)$, where $r=|p|$.
In particular, $j_H(t)=0$.
Set
\[
z=\sqrt{1+r^2},
\qquad
\phi(t)=\Phi_H(t).
\]
Then $h$ satisfies
\begin{equation}\label{radialfrozenPDE}
\partial_t h
=
\phi(t)z\partial_r^2h+
\left[
\phi(t)\left(\frac{2z}{r}+\frac r z\right)
+\frac{\mu r}{z}
\right]\partial_r h+
\mu\left(\frac2z+\frac1{z^3}\right)h.
\end{equation}
Set $q=\partial_r h$. Differentiating
\eqref{radialfrozenPDE} gives
\begin{equation}\label{radialderivativePDE}
\begin{aligned}
\partial_t q
={}&
\phi z\partial_r^2q
+
\left[
\phi\left(\frac{2z}{r}+\frac{2r}{z}\right)
+\frac{\mu r}{z}
\right]\partial_r q\\
&+
\left[
\phi\left(\frac1{z^3}-\frac2{zr^2}\right)
+\mu\left(\frac2z+\frac2{z^3}\right)
\right]q-
\mu r\left(\frac2{z^3}+\frac3{z^5}\right)h.
\end{aligned}
\end{equation}
Fix $T>0$. By Theorem~\ref{mainthm},
\[
0<\delta_T\leq\phi(t)\leq\frac{\mu}{3},
\qquad 0\leq t\leq T.
\]
Moreover, Lemmas~\ref{frozenestimates}
and~\ref{momentLipschitz}, applied to $H$
with frozen coefficients $(\phi,0)$, imply
that $\phi$ is H\"older continuous on $[0,T]$.
Throughout the following approximation argument,
we keep this coefficient function $\phi$ fixed.

We first prove preservation of radial monotonicity
for the frozen equation with smooth, compactly
supported initial data.
Let $f\in C_c^\infty(\R^3)$ be a nonnegative,
radially nonincreasing probability density,
and let $u$ be the solution of
\[
\partial_tu
=
\phi(t)\nabla_p\cdot(D\nabla_pu)
+
\mu\nabla_p\cdot(vu),
\qquad
u(0,p)=f(p).
\]
For $R$ sufficiently large that
$\operatorname{supp}f\subset B_R$, let $u^R$
solve the corresponding Dirichlet problem
\[
\begin{cases}
\displaystyle
\partial_tu^R
=
\phi(t)\nabla_p\cdot(D\nabla_pu^R)
+
\mu\nabla_p\cdot(vu^R),
& (t,p)\in(0,T]\times B_R,\\[1mm]
u^R(t,p)=0,
& (t,p)\in[0,T]\times\partial B_R,\\
u^R(0,p)=f(p),
& p\in B_R.
\end{cases}
\]
Classical linear parabolic theory gives a
nonnegative solution whose spatial gradient
is continuous up to the initial time and
the boundary. Here the compatibility conditions
hold because $f$ vanishes near $\partial B_R$.
Rotational invariance and uniqueness give
\[
u^R(t,p)=k^R(t,|p|).
\]
Set
\[
q^R(t,r)=\partial_r k^R(t,r).
\]
The function $q^R$ satisfies
\eqref{radialderivativePDE}, with $h$ replaced
by $k^R$ and $q$ replaced by $q^R$.
Its forcing term is nonpositive because
$k^R\geq0$.
The coefficient multiplying $q^R$ satisfies
\[
\phi(t)\left(\frac1{z^3}-\frac2{zr^2}\right)
+
\mu\left(\frac2z+\frac2{z^3}\right)
\leq
\phi(t)+4\mu
\leq\frac{13\mu}{3}.
\]
Smoothness at the center gives
\[
q^R(t,0)=0.
\]
At the outer boundary, nonnegativity of $k^R$
and $k^R(t,R)=0$ imply
\[
q^R(t,R)\leq0.
\]
Moreover,
\[
q^R(0,r)\leq0
\]
by the radial monotonicity of $f$.
Apply the maximum principle to
$e^{-5\mu t}q^R$.
Its zeroth-order coefficient is strictly negative,
and its forcing term is nonpositive.
A positive maximum cannot occur at $r=0$,
at $r=R$, or at the initial time.
At any $r>0$,
the coefficients of
\eqref{radialderivativePDE} are regular,
and the equation gives a contradiction.
Consequently,
\[
\partial_r k^R(t,r)\leq0,
\qquad
0\leq t\leq T,\quad 0<r<R.
\]
Extending $u^R$ by zero outside $B_R$ therefore
gives a radially nonincreasing function on $\R^3$.

We next pass to the whole space.
Let $P_t$ solve the frozen stochastic differential
equation with coefficients $(\phi,0)$ and initial
distribution $f(p)\,dp$, and define
\[
\tau_R=\inf\{s\geq0:\ |P_s|\geq R\}.
\]
The Dirichlet solution is the density of this
diffusion killed upon leaving $B_R$:
\[
\int_Au^R(t,p)\,dp
=
\mathbb P(P_t\in A,\ t<\tau_R)
\]
for every Borel set $A\subset\R^3$,
where $u^R$ denotes its zero extension.
Since $u(t,\cdot)$ is the density of $P_t$,
we have $0\leq u^R(t,\cdot)\leq u(t,\cdot)$
almost everywhere, and hence
\[
\|u(t)-u^R(t)\|_{L^1}
=
1-\mathbb P(t<\tau_R)
=
\mathbb P(\tau_R\leq t).
\]
The diffusion is global and nonexplosive.
Its continuous paths therefore imply
$\tau_R\to\infty$ almost surely as $R\to\infty$.
Consequently,
\[
\sup_{0\leq t\leq T}
\|u(t)-u^R(t)\|_{L^1}
\leq
\mathbb P(\tau_R\leq T)
\longrightarrow0.
\]
The cone of nonnegative, radially nonincreasing
functions is closed in $L^1(\R^3)$.
Thus $u(t,\cdot)$ is radially nonincreasing
for every $t\in[0,T]$.

Finally, choose nonnegative, radially
nonincreasing probability densities
$f_n\in C_c^\infty(\R^3)$ such that
\[
f_n\longrightarrow H_0
\qquad\text{in }L^1(\R^3).
\]
Such an approximation is obtained by radial
truncation, convolution with a smooth radially
nonincreasing mollifier, and normalization.
Let $u_n$ solve the frozen equation with
coefficients $(\phi,0)$ and initial data $f_n$.
The preceding argument shows that
$u_n(t,\cdot)$ is radially nonincreasing.
Since $j_H=0$, the given nonlinear solution $H$
also solves this same frozen equation, with
initial data $H_0$.
The contraction estimate
\eqref{frozenL1contraction} therefore gives
\[
\|u_n(t)-H(t)\|_{L^1}
\leq
\|f_n-H_0\|_{L^1},
\qquad 0\leq t\leq T.
\]
Passing to the limit shows that $H(t,\cdot)$
is radially nonincreasing.
By Theorem~\ref{regularizationthm}, $H$ is smooth
for positive times, so this implies
\[
\partial_r h(t,r)\leq0,
\qquad 0<t\leq T,\quad r>0.
\]
Since $T>0$ was arbitrary, the conclusion follows.
\end{proof}

We now establish the comparison principle
for two even solutions.

\begin{lemma}\label{evencomparisonlemma}
Let $G$ and $H$ be global weak solutions
corresponding to even probability densities
$G_0$ and $H_0$ with finite first moments.
Assume that $H_0$ is radially nonincreasing
and that
\begin{equation}\label{eveninitialorder}
F_G(0,z)\geq F_H(0,z),
\qquad z\geq1.
\end{equation}
Then
\begin{equation}\label{evencomparisonorder}
F_G(t,z)\geq F_H(t,z),
\qquad t\geq0,\ z\geq1.
\end{equation}
Consequently,
\begin{equation}\label{evenPhicomparison}
\Phi_G(t)\geq\Phi_H(t),
\qquad t\geq0.
\end{equation}
\end{lemma}

\begin{proof}
Set
\[
U(t,z)=F_G(t,z)-F_H(t,z).
\]
Recall the definition~\eqref{energyLnew} of the operator $\mathcal{L}_\phi$. 
Subtracting the equations~\eqref{energyFnew}
for $F_G$ and $F_H$, we obtain
\begin{equation}\label{evenUequation}
\partial_tU
=
\mathcal L_{\Phi_G(t)}U
+
(\Phi_G(t)-\Phi_H(t))K_H,
\end{equation}
where
\[
K_H(t,z)
=
(z-z^{-1})\partial_z^2F_H
-
(2-z^{-2})\partial_zF_H.
\]
By Lemma~\ref{radialmonotonicitylemma},
\[
H(t,p)=h(t,|p|),
\qquad
\partial_rh(t,r)\leq0.
\]
Writing $r=\sqrt{z^2-1}$, we have
\[
\rho_H(t,z)=4\pi zr\,h(t,r),
\qquad
\partial_zF_H=-\rho_H.
\]
A direct calculation gives
\[
K_H(t,z)
=
-4\pi zr^2\partial_rh(t,r)
\geq0.
\]
Moreover,
\[
\int_1^\infty K_H(t,z)\,dz=
-4\pi\int_0^\infty
r^3\partial_rh(t,r)\,dr=
12\pi\int_0^\infty r^2h(t,r)\,dr
=3.
\]
The boundary terms vanish because
$h$ is nonnegative, radially nonincreasing,
and integrable.
Define
\[
U_-(t,z)=\max\{-U(t,z),0\},
\qquad
A(t)=\int_1^\infty z^{-3}U_-(t,z)\,dz.
\]
By~\eqref{PhiFcumulative},
\[
\Phi_G(t)-\Phi_H(t)
=
\frac{2\mu}{3}
\int_1^\infty z^{-3}U(t,z)\,dz,
\]
and therefore
\[
(\Phi_H(t)-\Phi_G(t))_+
\leq\frac{2\mu}{3}A(t).
\]
We next prove the inequality
\begin{equation}\label{evenKato}
\partial_tU_-
\leq
\mathcal L_{\Phi_G(t)}U_-
+
\frac{2\mu}{3}A(t)K_H(t,z)
\end{equation}
in the distributional sense on
$(0,\infty)\times(1,\infty)$.
For $\varepsilon>0$, set
\[
\beta_\varepsilon(s)
=
\frac{\sqrt{s^2+\varepsilon^2}-s}{2}.
\]
Then $\beta_\varepsilon\to(\,\cdot\,)_-$,
$-1\leq\beta_\varepsilon'\leq0$, and
$\beta_\varepsilon''\geq0$.
By
\eqref{evenUequation}, and using $K_H\geq0$, 
\[
\begin{aligned}
(\partial_t-\mathcal L_{\Phi_G(t)})
\beta_\varepsilon(U)
&=
\beta_\varepsilon'(U)(\Phi_G-\Phi_H)K_H
-
\Phi_G(t)(z-z^{-1})
\beta_\varepsilon''(U)|\partial_zU|^2\\
&\leq
(\Phi_H-\Phi_G)_+K_H
\leq
\frac{2\mu}{3}A(t)K_H.
\end{aligned}
\]
Letting $\varepsilon\downarrow0$ in the sense
of distributions on
$(0,\infty)\times(1,\infty)$
yields \eqref{evenKato}.

Let $\mathcal L_\phi^*$ denote the formal
adjoint of $\mathcal L_\phi$ with respect
to Lebesgue measure. A direct calculation
gives
\[
\mathcal L_\phi^*(z^{-3})
=
\frac{3\mu z^2-5\mu-15\phi}{z^6}
\leq3\mu z^{-3}.
\]
For positive times, smoothness of
$G$ and $H$ implies
\[
U_-(t,z)=O((z-1)^{3/2})
\qquad\text{as }z\downarrow1,
\]
locally uniformly in time.
Together with $0\leq U_-\leq1$,
this justifies integration of
\eqref{evenKato} against $z^{-3}$
by standard localization.
Thus, for $0<s<t$,
\[
A(t)
\leq A(s)
+3\mu\int_s^t A(\tau)\,d\tau+
\frac{2\mu}{3}
\int_s^t A(\tau)
\left(
\int_1^\infty z^{-3}
K_H(\tau,z)\,dz
\right)d\tau.
\]
Using $K_H\geq0$ and
\[
\int_1^\infty K_H(\tau,z)\,dz=3,
\]
we obtain
\begin{equation}\label{evenGronwall}
A(t)\leq A(s)+5\mu\int_s^t A(\tau)\,d\tau.
\end{equation}

Finally, $A$ is continuous on $[0,\infty)$.
Indeed,
\[
|A(t)-A(s)|
\leq\frac12\bigl(
\|G(t)-G(s)\|_{L^1}+\|H(t)-H(s)\|_{L^1}
\bigr).
\]
By~\eqref{eveninitialorder}, $A(0)=0$.
Letting $s\downarrow0$ in
\eqref{evenGronwall} and applying
Gronwall's inequality yields
\[
A(t)=0,
\qquad t\geq0.
\]
Hence $U_-(t,z)=0$ almost everywhere
in $z$, for every $t\geq0$.
By continuity of the tail distributions,
\[
F_G(t,z)\geq F_H(t,z),
\qquad t\geq0,\ z\geq1.
\]
Finally,~\eqref{PhiFcumulative} gives
\[
\Phi_G(t)-\Phi_H(t)
=
\frac{2\mu}{3}
\int_1^\infty z^{-3}
\bigl(F_G(t,z)-F_H(t,z)\bigr)\,dz
\geq0.
\]
This completes the proof.
\end{proof}

We can now complete the proof of the theorem.

\begin{proof}[Proof of Theorem~\ref{uniformPhitheorem}]
The conclusion for $1<E<2$ is proved in Lemma~\ref{lowenergyPDElemma}. 
We now establish uniform positivity of
$\Phi_G$ for arbitrary finite energy $E>1$.
Let $G_0^*$ be the symmetric decreasing
rearrangement of $G_0$.
Then $G_0^*$ is a radially nonincreasing
probability density satisfying
\begin{equation}\label{rearrangementmass}
\int_{|p|\leq r}G_0^*(p)\,dp
\geq
\int_{|p|\leq r}G_0(p)\,dp,
\qquad r>0.
\end{equation}
Since $\langle p\rangle$ is radially
increasing, it follows that
\begin{equation}\label{rearrangementenergy}
\int_{\R^3}\langle p\rangle G_0^*(p)\,dp
\leq E.
\end{equation}
For $0<\varepsilon<1$, define
\[
H_0(p)=\varepsilon^{-3}G_0^*(p/\varepsilon).
\]
Then $H_0$ is a radially nonincreasing
probability density with finite first moment.
By~\eqref{rearrangementmass},
\[
\int_{|p|\leq r}H_0(p)\,dp=
\int_{|p|\leq r/\varepsilon}G_0^*(p)\,dp\geq
\int_{|p|\leq r}G_0(p)\,dp.
\]
Equivalently,
\begin{equation}\label{rearrangementorder}
F_G(0,z)\geq F_H(0,z),
\qquad z\geq1.
\end{equation}
Let
\[
E_H=\int_{\R^3}\langle p\rangle H_0(p)\,dp.
\]
Since $\langle\varepsilon p\rangle
\leq
1+\varepsilon(\langle p\rangle-1)$,
we have
\[
E_H
=
\int_{\R^3}
\langle\varepsilon p\rangle G_0^*(p)\,dp\leq
1+\varepsilon
\int_{\R^3}
(\langle p\rangle-1)G_0^*(p)\,dp\leq1+\varepsilon(E-1).
\]
Since $H_0$ is a probability density,
we have $E_H>1$. Choosing
$\varepsilon>0$ sufficiently small,
we obtain $1<E_H<2$.
Let $H$ be the global weak solution
with initial data $H_0$.
By Lemma~\ref{lowenergyPDElemma},
there exists $\delta_H>0$ such that
\begin{equation}\label{referencePhigap}
\Phi_H(t)\geq\delta_H,
\qquad t\geq0.
\end{equation}
Moreover, Lemma~\ref{radialmonotonicitylemma}
shows that $H(t,\cdot)$ remains
radially nonincreasing.
Therefore the assumptions of
Lemma~\ref{evencomparisonlemma}
are satisfied, and
\[
F_G(t,z)\geq F_H(t,z),
\qquad t\geq0,\ z\geq1.
\]
Consequently,
\[
\Phi_G(t)\geq\Phi_H(t),
\qquad t\geq0.
\]
Combining this inequality with
\eqref{referencePhigap}, we conclude that
\[
\Phi_G(t)\geq\delta_H>0,
\qquad t\geq0.
\]
The constant $\delta_H$ depends only
on $H_0$ and $\mu$, and hence only
on $G_0$ and $\mu$.
This proves~\eqref{uniformPhibound}
and completes the proof.
\end{proof}

%===========================================================

\end{document}